\documentclass[10pt,english]{smfart}
\usepackage{times}

  \usepackage{epsfig}
\usepackage{amssymb}
\usepackage{amsmath}
\usepackage{amsthm}
\usepackage{graphics}
\usepackage{hyperref}

\usepackage{xcolor}

\newtheorem{theorem}{Theorem}
\newtheorem{conj}{Conjecture}
\newtheorem{question}{Question}
\newtheorem{problem}{Problem}
\newtheorem{example}{Example}
\newtheorem{definition}{Definition}
\newtheorem{remark}{Remark}

\newcommand{\tH}{\widetilde{H}}
\newcommand*\sq{\mathbin{\vcenter{\hbox{\rule{.3ex}{.3ex}}}}}

\begin{document}

\title{H\'enon maps: a list of open problems}
\author{Julia Xénelkis de Hénon}
\maketitle

\begin{abstract}
	We propose a set of questions on the dynamics of Hénon maps from the real, complex, algebraic and arithmetic points of view. 
\end{abstract}

\tableofcontents
\section{Introduction (C. Favre, T. Firsova, L. Palmisano, J. Raissy,  and G. Vigny (Eds.))}
A workshop `Dynamics of H\'enon maps: Real, Complex and Beyond" took place at BIRS, Banff
in April 2023. The purpose of this meeting was to bring together international experts working on various aspects of H\'enon maps. 
Recall that these maps are two-dimensional discrete dynamical systems
which are  ubiquitous in low dimensional dynamics, and among the most studied examples 
exhibiting chaotic behavior. 
Quadratic H\'enon maps 
\[ H_{a,c} (x,y) := (ay+x^2+c, a x) \]
are arguably the simplest examples.
Here $a$ and $c$ are fixed parameters and $x,y$ are affine coordinates\footnote{The exact definition of (quadratic) H\'enon maps may differ from section to section in these notes, as one might want to
conjugate them by affine transformations to exploit various aspects of the original equations.}. 
These maps can be analyzed  over the real numbers using techniques from smooth dynamical systems, or over the complex numbers
and then complex analysis and geometry play crucial roles.
They are amenable to generalizations, by replacing $x^2+c$ by higher degree polynomials, or even transcendental maps, and we may consider 
finite composition of such maps. 
We may also consider them with coefficients in number fields, and look at them from the perspective of arithmetic dynamics. 
Many of the most recent breakthroughs were actually made by combining several techniques coming from these different fields. It was delightful to attend series of 
talks blending so many different ideas. Many interesting questions were raised during the conference,  a fact which encouraged us to
collect them in a single text. 

\medskip

H\'enon introduced his family of maps in the real domain as a simplified model of the Poincar\'e section of
the first return map of the Lorenz flow~\cite{henon}. 
In 1976, H\'enon made numerical experiments for the map $H_{\sqrt{0.3}, -1.4}$\footnote{H\'enon actually considered the map $h(x,y)=(1-1.4\cdot  x^2+y\, , \, 0.3\cdot  x)$ which is affinely conjugate to it.}
and observed
that an initial point of the plane either approaches
a set of points known since then as the H\'enon strange attractor, or diverges to infinity under iterations. The H\'enon
attractor has a fractal nature: it is smooth in the unstable direction and has a Cantor-like structure in the transversal direction. This led H\'enon 
to conjecture the existence of an ergodic measure which restricts to the Lebesgue measure in the smooth direction (a.k.a.
an  SRB measure). In 1981, Jakobson \cite{MR0630331} proved the existence of a set of positive Lebesgue measure of parameters $c$ for which $x^2+c$ displays an SRB measure. In the 90's, Benedicks and Carleson \cite{BC} reworked Jakobson's theorem and further generalized it to describe the dynamics of H\'enon maps. They proved $H_{a,c}$ display a strange attractor for $a$ small and for a set of parameters $c$ of positive Lebesgue measure.  Benedicks-Carleson's breakthrough has been further developed by Mora-Viana \cite{MoraViana}, Wang-Young \cite{WY}, and Takahasi \cite{Ta}.

In 1996, during his inaugural lecture at Coll\`ege de France, Yoccoz proposed an alternative approach to prove the H\'enon conjecture, with Sinai's positive entropy conjecture lying in the horizon.  To this end, he introduced a combinatorial and topological approach, based on the notion of strong regularity, that he used to give yet another proof of Jakobson's theorem \cite{Yo19}. The second step of Yoccoz' program was completed more recently by Berger in \cite{Berger} who generalized this notion of strong regularity,  leading him in particular to an alternative proof of Benedicks-Carleson theorem.

\medskip

The theory of H\'enon maps in the complex domain started with the seminal work of 
 Friedland and Milnor \cite{FM}, who used Jung's theorem to show that 
every polynomial automorphism of the complex affine plane is affinely conjugate to either an affine map,
or a map preserving the pencil $x=\mathrm{cst}$, or to a finite compositions of generalized H\'enon maps (the latter class being usually called
complex H\'enon maps nowadays). 
In the early 90's Hubbard and his collaborators developed a topological approach to the study of Henon maps giving description of the Fatou sets \cite{HOV,HPV}. Hubbard pointed out that H\'enon maps appeared as natural generalizations of quadratic polynomials (by taking $a \to  0$ in $H_{a,c}$), which he used with Oberste-Vorth to give topological description of Julia sets~\cite{MR1351520}.

An important breakthrough was the introduction of pluripotential techniques to construct invariant currents by Forn{\ae}ss and Sibony~\cite{FS} and by Bedford and Smillie~\cite{BS}. 
The latter authors, partly with Lyubich, further developed in a series of influential papers (e.g.,~\cite{BSIII,BLS,BLS2})  a thorough study of the ergodic properties
of H\'enon mappings and of their stability properties linked to hyperbolicity.
These works were pursued and generalized to other invertible maps by many others including  Diller, Cantat, Dujardin, Dinh and Sibony~\cite{Diller96, Cantat,Dujardin_duke,MR2045508,Dinh_decay, Dinh_Sibony_horizontal}. Very recently, Bianchi and Dinh
\cite{BianchiDinh1} made significant progress in the study of the fine statistical properties of the maximal entropy measure.

It is intriguing to see in retrospect how these seemingly simple maps have produced such an elaborate and successful theory. Note however that
the results obtained so far are most complete in the case of dissipative maps (that is $|a|<1$), while the understanding of the conservative case (that is $|a|=1$) remains less developed.

\medskip

In the past decade, the algebraic and arithmetic aspects of dynamical systems defined by rational maps
have  also been developed extensively. We refer to the survey~\cite{survey_arithmetic_dynamics} in which one can find a large set of open problems in this
emerging field. For H\'enon maps defined over a number field or over a function field,  Silverman \cite{S}, and later Kawaguchi \cite{K},
constructed a suitable height function that lead to interesting analogs
of the Northcott property (see also \cite{I}). H\'enon maps have also been the testing ground of some  important conjectures in arithmetic dynamics
like the Kawaguchi-Silverman conjecture on arithmetic degrees \cite{Kawaguchi_silverman}, 
or the dynamical Manin-Mumford problem which was partially solved by Dujardin and Favre
using both height and Pesin theories \cite{DF2}. Deep connections exist between the arithmetic of these systems and pluripotential theoretic techniques:
it is for instance possible to retrieve the equidistribution of repelling periodic orbits using a theorem on the equidistribution of points of small height by Yuan \cite{Yuan}.

\medskip

The study of H\'enon maps is still very active, as shown by the recent breakthroughs 
in the study of wandering domains. On the one hand, Ou~\cite{Ou}
has proved the absence of wandering domains for strongly dissipative doubly infinitely period-doubling renormalizable real H\'enon maps. 
On the other hand, 
Berger and Biebler \cite{BB}  exhibited wandering domains for complex
H\'enon maps in 2023 by mixing deep techniques coming  from both
real and complex dynamics. We also witness exciting new developments extending the already rich theory
of H\'enon maps to more general systems such as  transcendental diffeomorphisms of the complex plane~\cite{ABFH, ABFH3, ABFH2}, or
higher dimensional invertible rational maps, \cite{Dethelin_Vigny,Dinh_Sibony_JFA, Dinh_Sibony_Kahler, Dinh_Sibony_equi,Guedj_Sibony},   
where questions arising from complex dynamics led to profound developments in complex geometry such as PB currents, density currents and superpotential theory. Several developments have been made also in the context of higher dimensional unfoldings of homoclinic tangencies with a rank one saddle point, \cite{PV,V,WY1}.
We hope that gathering these questions and open problems at one place will reinforce the community
 and  attract new generations of researchers to work on these beautiful and rich 
objects.   

\paragraph*{Acknowledgments} The idea to collect this list of problem arose after holding a problem session at a BIRS conference which was lively chaired by M. Abate.
We thank all participants of this session and the other contributors to this list to have generously shared their ideas and problems
on Hénon maps. We are also grateful to Zin Arai, André de Carvalho, Sébastien Gouëzel, Stéphane Lamy, Joseph Silverman, and Liz Vivas for their careful reading and their comments on previous versions of this manuscript.



\section{Real H\'enon maps (P. Berger)}

We propose a set of questions on the dynamics of H\'enon maps in the real domain, or more generally
on entire diffeomorphisms of $\mathbb{R}^2$.

\subsection{Strange attractors}
Attractors play an important role in the study of dynamical systems since the 60's
(Lorenz attractor, H\'enon strange attractor, etc.). This notion is quite flexible and
can cover many different situations in which a substantial set of points (either in the topological
or measurable sense) is converging to some invariant compact subset. We refer to Milnor~\cite{MR0790735}
for a discussion of various possible definitions of attractors.  For instance,  a measure-theoretical attractor is an invariant compact subset which attracts a set of Lebesgue positive measure  and which is minimal with this  property.

In the case of unimodal interval dynamics, measure-theoretical attractors can be classified into four types: cyclic, solenoidal, interval and wild (see, e.g.,~\cite{MR1132757}).
The  latter two classes are arguably the most interesting.
Jakobson~\cite{MR0630331} proved the abundance of quadratic maps displaying a stochastic interval of attractors (induced by an SRB measure). On the other hand,  Bruin-Keller-Nowicki-van Strien~\cite{MR1370759} showed the existence of a polynomial unimodal map 
displaying a wild attractor:  an invariant Cantor set attracting Lebesgue almost every point and included in a transitive interval.

Van Strien~\cite[Question~1.9]{MR2600680} asked whether a suitable analog of wild attractors could exist
for H\'enon maps (of some degree). More precisely, we ask:
\begin{question} \normalfont
Can we find a H\'enon map (of some degree) which admits
a wild attractor, i.e., a Cantor set which attracts a set of Lebesgue positive measure and which is strictly included in a transitive set ?
\end{question}

Returning to stochastic attractors, observe that
the existence of a parameter $c\in \mathbb{R}$ for which the quadratic map $x \mapsto x^2+c$
displays an absolutely continuous measure is easy to ensure. It suffices to pick
a parameter $c$ such that the post-critical orbit is finite but not periodic.
In dimension $2$, the existence of a positive measure set of parameters of H\'enon maps displaying an attractor supporting an invariant SRB measure
\footnote{i.e., a measure whose conditional measures along unstable curves are absolutely continuous.}
is a fundamental result, whose proof still remain difficult and lengthy (\cite{BC,Berger,Ta,WY}).

\begin{question} \normalfont

    Is there a quick proof for the existence of SRB for some parameters of the H\'enon maps?

\end{question}

A positive answer to this question might help finding new examples of stochastic attractors in the H\'enon maps.

\subsection{Non-statistical behavior}
Let $f$ be any smooth dynamical system. 
We say that a point $x$ has a non-statistical  behavior (or simply is non-statistical)
if its sequence of empirical measures  $e_n(x):= \tfrac1n \sum_{k=0}^{n-1} \delta_{x_n} $
is not converging. We say that $f$ is non-statistical if there is a
positive Lebesgue measure set of points with non convergent empirical measures. 
Ruelle~\cite{MR1858471} asked whether non-statistical dynamics could exist
persistently\footnote{Ruelle used another terminology and talked about historical behavior.}.

In polynomial  dynamics, two phenomena give rise to non-statistical dynamics. 
The first one was discovered by  Hofbauer and Keller: 

\begin{theorem}[\cite{HK90}]\label{thm:hofbauer}
There exist uncountably many $c\in \mathbb{R}$ such that the quadratic polynomial  $P_c(x):= x^2+c$ has non-statistical dynamics. 
More precisely, Lebesgue almost every non-escaping point  $x$ has non-statistical  behavior.
\end{theorem} 
In \cite{Ta22},  Talebi gave a  counterpart of this result for rational functions on the Riemann sphere. 
In these two results, the set of non-statistical points is of full measure, but of empty interior.

The  second occurrence of non-statistical dynamics is related to the notion of wandering stable component that we now introduce.
\begin{definition}
	A {stable domain} of $f$ is a connected open subset $U$ such that \[\lim_{n\to\infty}d(f^n(x), f^n(y))= 0\] for all $x,y\in U$. A {stable component} is a maximal stable domain.
	A stable component is {wandering} if it is not preperiodic. 
\end{definition} 
Colli and Vargas \cite{CV01}
gave the first example of a smooth dynamical system  having a wandering stable component formed by points with non-statistical behavior. In \cite{KS17},  Kiriki and Soma constructed a locally dense set of such dynamics in the $C^r$-category with $r<\infty$.  This also occurs for polynomial maps: 
\begin{theorem}[\cite{BB}]
	There is a locally dense set of real sextic polynomials $P(x)=  x^6+ a_4 x^4 + \cdots + a_0$, 
 for which the map $(x, y)\mapsto (P(x)-b y,x)$ displays a wandering stable component containing only points with non-statistical  behavior.  
\end{theorem} 

The proof of this theorem actually implies the existence of a wandering Fatou component at the same parameters for its complex counterpart.
This is in sharp contrast with the one-dimensional situation for which no wandering Fatou components exist by Sullivan's theorem. 
 
 Conversely, we can ask  whether a counterpart of Hofbauer-Keller phenomenon appears within the H\'enon family. 
We can formulate this question in more precise terms. 
\begin{question}  \normalfont
	Do there exist  $\varepsilon>0$ and a locally dense set\footnote{i.e., whose closure has non-empty interior.} $E$ of quadratic H\'enon maps for which 
	every point starting in some set of Lebesgue measure at least $\varepsilon$ has non-statistical behavior? 
\end{question} 

By \cite[Theorem 1.14]{Ta22}, this would imply the existence of a generic set of the closure of $E$  of non-statistical dynamics.  

\subsection{Conservative dynamics}
Let $f\colon S \to S$ be any homeomorphism of a closed surface $S$. 
An annular rotation domain for $f$ is by definition 
an $f$-invariant open annulus that does not contain any periodic points.
Such domains play an important role in conservative dynamics.
The next result can be deduced from the works of \cite{FL03, K10, KLN15, X06} (Le Calvez, private communication). 
\begin{theorem}
	Let $S$ be a closed surface  endowed with a symplectic form.   Let $f\colon S \to S$ be any 
symplectomorphism 	of class $C^r$ with $r\ge 1$ that contains at least one periodic point and satisfying the following conditions.
	\begin{enumerate}
		\item Any eigenvalue of any periodic point does not belong to $\{e^{2\pi ip/q}: p/q\in \mathbb Q\}$.
		\item For every hyperbolic periodic points $P, Q\in \mathrm{Per}(f)$, $W^s(P)$  is transverse to $W^u(Q)$.
		\item Every elliptic point $P\in \mathrm{Per}(f)$ is surrounded by arbitrarily close KAM circles.   
		\item There are no annular rotation domain.
	\end{enumerate}
	Then  $\bigcup_{P\in \mathrm{Per}(f)} W^s(P)$  is dense in $S$.  
\end{theorem}

It is natural to ask whether the third condition is superfluous. The recent result \cite{OC23} suggests that this may be the case.

\begin{question}\normalfont
	Is there a real, conservative, polynomial H\'enon map with an annular rotation domain?
	Is there an open set of such real, conservative, polynomial H\'enon maps? 
\end{question} 

An annular rotation domain is said to be trivial when the whole dynamics is conjugate to a rotation.

\begin{question}\label{entireR2}\normalfont
	Is there an entire\footnote{i.e., an analytic map which extends to a holomorphic map of $\mathbb{C}^2$.} symplectomorphism of  $\mathbb{R}^2$ with a non trivial annular rotation domain?  
\end{question} 
Recently, an entire map of the cylinder $\mathbb{R}\times \mathbb{R}/\mathbb{Z}$ having a bounded rotation domain on which the dynamics is not conjugate to a rotation
has been exhibited in  \cite{Berg22}, thereby disproving a conjecture by  Birkhoff \cite{Bi41} .
The construction in \cite{Berg22} also gives an example of a symplectic entire automorphism of $\mathbb{C}\times \mathbb{C}/\mathbb{Z}$ without periodic point and with a non-empty set of (Lyapunov) unstable points\footnote{an orbit $(x_n)_{n\ge0} $ is Lyapunov stable  if  for any $y_0$ close enough to $x_0$, then $y_n$ stays close to $x_n$ for all $n\ge0$.}. 
 Hence Question \ref{entireR2}  may shed light on the following intriguing problem.
\begin{problem} \normalfont
	Construct an entire symplectic automorphism of $\mathbb{C}^2$ without periodic point and with non-empty set of (Lyapunov) unstable points with bounded orbit. 
\end{problem}

\section{Dissipative real H\'enon maps (S. Crovisier and E. Pujals)}

\subsection{Mild dissipation}
Quadratic real H\'enon maps
\[f_{c, b} (x, y):= (x^2+c-by, x),\]
with Jacobian $b$ close to zero share some properties of the quadratic family on the interval: some results are obtained by perturbative methods (for instance~\cite{BC,dCLM}) and are known when $|b|$ is extremely tiny. In~\cite{CP} another approach has been introduced which allows to reduce the dynamics to a one-dimensional system.

\begin{definition}
The map $f_{c,b}$ is \emph{mildly dissipative} if it is dissipative (i.e., $|b|<1$) and if for any ergodic invariant measure $\mu$
which is not supported on a sink, and for $\mu$-almost every point $x$,
both components of the stable curve $W^s(x)\setminus \{x\}$ are unbounded.
\end{definition}

Under this assumption, and in restriction to any open topological disc $D\subset \mathbb{R}^2$ that is compactly mapped inside itself, the dynamics is semi-conjugated to a (non-trivial) continuous map of a real tree. Other strong consequences can be derived (e.g., a closing lemma, or a description of zero entropy dynamics, see below). Using Wiman's theorem (in the same spirit like in~\cite{DL, MR3213832}), one can prove that H\'enon maps
are mildly dissipative once $|b|<1/4$. In this way we obtain dynamical informations for \emph{all} H\'enon maps
having their Jacobian in $(-1/4,1/4)$ and not only for those satisfying $|b|\ll 1$.
One expects that this property extends beyond the bound obtained through Wiman's theorem.

\begin{question}
Which real H\'enon maps $f_{c,b}$ are mildly dissipative?
Is this property satisfied by all H\'enon maps with $|b|<1$?
\end{question}

In some cases \cite{CP} proves  that the mild dissipation is an open property,
but we don't know if this holds in general.

\subsection{Maps with zero entropy}
It is well-known that the quadratic maps $f_c(x):=x^2+c$ have their topological entropy equal to zero exactly when $c$ belongs to some interval
$(-\infty,c_0]$. At the critical parameter $c_0$, the dynamics exhibits an odometer, which is the limit set of an infinite sequence of successive renormalizations of period $2$. This result persists inside any line $c\mapsto f_{c,b}$, provided $|b|$ is smaller than some number $\varepsilon>0$ small, as it has been shown in~\cite{dCLM}.
Let us consider the locus in the parameter space where the topological entropy vanishes:
\[\mathcal{E}_0:=\{(b,c)\in \mathbb{R}^2, h_{\text{top}}(f_{c,b})=0\}.\]
Inside the strip $(-\varepsilon,\varepsilon)\times \mathbb{R}$ this set is bounded by an analytic arc $\{(b,c_0(b)),|b|<\varepsilon\}$.
Moreover any map $f_{c,b}$ with $c<c_0(b)$ can be renormalized at most finitely many times;
and for an open and dense subset of these parameters, the dynamics is Morse Smale
(i.e., is supported by finitely many hyperbolic periodic orbits).
When $c=c_0(b)$, the sequence of renormalizations converges towards a particular unimodal map of the interval.

When $|b|$ is larger but smaller than $1/4$ (so that it is mildly dissipative),
\cite{CPT} describes the dynamics of $f_{c,b}\in \mathcal{E}_0$.
In particular, all maps $f_{c,b}\in \partial \mathcal{E}_0$ are infinitely renormalizable (with renormalization periods eventually equal to $2$), solving a conjecture by Tresser (which is still open when we don't assume the mild dissipation).
One may wonder if the converse holds.

\begin{question}
Let us consider any infinitely renormalizable mildly dissipative map $f_{c,b}\in\mathcal{E}_0$.
Is it the limit of maps with positive entropy? Is it the limit of Morse-Smale maps?
\end{question}

It is also natural to try to describe the boundary $\partial\mathcal{E}_0$: is it a (piecewise smooth) arc?
One would like to implement the strategy developed in~\cite{dCLM} for $b$ close to $0$:

\begin{question}
Let us consider any infinitely renormalizable mildly dissipative map in $\mathcal{E}_0$.
Does the sequence of renormalizations converge?
\end{question}

We also don't know if different combinatorics of the renormalizations may occur.

\begin{question}
Does there exist a mildly dissipative map $f_{c,b}\in \mathcal{E}_0$
which admits arbitrarily deep renormalizations with odd periods?
\end{question}

\subsection{Set of periodic points}
One would like to describe the dynamics through its periodic orbits.
To any periodic point $p$, one associates two Lyapunov exponents
$\lambda^-(p)\leq \lambda^+(p)$. When they do not vanish and have different sign,
we say that $p$ is a saddle.

As mentioned previously, for mildly dissipative real H\'enon maps the set of periodic points is dense
in the union of the supports of the invariant probability measures~\cite{CP}.
The same property holds for any complex H\'enon maps~\cite{MR4142445}.

\begin{question}
For any real H\'enon map, does the closure of the set of periodic orbits support all the invariant probability measures?
\end{question}

The next step is to describe how periodic saddles are organized.
We say that two saddles $p,q$ are homoclinically related if
there exists $k\geq 0$ such that the invariant curves
$W^u(p), W^s(f^k(q))$ (and  $W^s(p), W^u(f^k(q))$  as well) intersect transversally. 
This defines an equivalence relation which decomposes the set of periodic saddles
into its homoclinic classes.
There may exist infinitely many periodic saddles which are not homoclinically related,
but we conjecture that their hyperbolicity should drop.

\begin{question}
For any map $f_{c,b}$ and any infinite set of periodic saddles $(p_n)$
which are pairwise not homoclinically related, do we have
$\min(|\lambda^-(p_n)|,|\lambda^+(p_n)|)\underset{n\to \infty}\longrightarrow 0$?
\end{question}

This questions goes beyond H\'enon maps, but \cite{CPT} implies a positive answer
in the particular case of mildly dissipative $f_{c,b}\in \mathcal{E}_0$.

\section{Symbolic dynamics for real H\'enon and Lozi maps (S. \v Stimac)}

Kneading theory is a combinatorial tool to understand the dynamics of a piecewise monotone map from the interval to itself and was developed by Milnor and Thurston~\cite{MilThu}. 
Applications extend from the topological classification to the computation of the entropy, to the counting of periodic orbits, and the construction of measures of maximal entropy. 
We propose several problems connected to the extension of this theory to
real H\'enon maps and the Lozi maps $\tH_{a,b}, L_{a,b} \colon \mathbb{R}^2 \to \mathbb{R}^2$, \[\tH_{a,b}(x,y) = (1 + y - ax^2, bx), \quad L_{a,b}(x,y) = (1 + y - a|x|, bx),\]
respectively\footnote{observe that $\tH_{a,b}$ is affinely conjugated to $H_{\sqrt{b},-a}$ from the introduction}. The Lozi maps are piecewise affine map that display the same  fold and bend behavior as the H\'enon maps, but are usually easier to analyze technically~\cite{Lozi1,Lozi2}.

\medskip

In \cite{MS}, the authors developed a kneading theory for the Lozi maps $L_{a,b}$ with $(a,b) \in \mathcal{M}$, where $\mathcal{M} = \{ (a,b) \in \mathbb{R}^2 : b > 0, \ a\sqrt{2} - b > 2, \ 2a + b < 4 \}$ is a set of parameters for which Misiurewicz proved the existence of a strange attractor (for details see \cite{MS} and \cite{MS2}). A kneading sequence $\bar k$ is defined as the itinerary of a turning point $T$, where turning points are points of transversal intersections of the $x$-axis and the unstable manifold $W^u$ of the fixed point $X$ of the attractor. Any kneading sequence is a bi-infinite sequence of $+$ and $-$.

The kneading set $\mathfrak{K} = \{ \bar k^n : n \in \mathbb{Z}\}$ is the set of all kneading sequences $\bar k^n$, $n \in \mathbb{Z}$, and every kneading sequence $\bar k = \bar k^n$, for some $n \in \mathbb{Z}$, has the following form: \[\bar k = \hspace{.2cm}+^{\hspace{-.5cm}\infty}\hspace{.2cm} w \pm \sq  \overrightarrow k_{\hspace{-.1cm} 0},\] where $\hspace{.2cm}+^{\hspace{-.5cm}\infty}\hspace{.2cm} = \dots + + +$, $w = w_0 \dots w_m$, for some $m \in \mathbb{N}_0$, $\overrightarrow k_{\hspace{-.1cm} 0} = k_0 k_1 k_2 \dots$, $w_0 = -$, $k_0 = +$, $w_i, k_j \in \{ -, + \}$ for $i = 1, \dots, m$ and $j \in \mathbb{N}$, and the little black square $\sq$ indicates where the 0th coordinate is located. Here for $\pm$ one can substitute any of $+$ and $-$.

In \cite{MS}, the authors prove that $\mathfrak{K}$ characterizes all itineraries of all points of the attractor of $L_{a,b}$. The proof is given in two steps. 
We say that an itinerary is $W^u$-admissible if it is realized 
by a point on the unstable manifold $W^u$. We first have:
\begin{theorem}
	A sequence $\hspace{.2cm}+^{\hspace{-.5cm}\infty}\hspace{.2cm} \overrightarrow p_{\hspace{-.1cm} n}$, where $\overrightarrow p_{\hspace{-.1cm} n} = p_n p_{n+1} \dots$ such that $p_n = -$ for some $n \in \mathbb{Z}$, is $W^u$-admissible if and only if for every kneading sequence $\hspace{.2cm}+^{\hspace{-.5cm}\infty}\hspace{.2cm} w \pm \sq\overrightarrow k_{\hspace{-.1cm} 0}$, such that $w = p_np_{n+1} \dots p_{n+m}$ for some $m \in \mathbb{N}_0$, we have $\sigma^{m+2}(\overrightarrow p_{\hspace{-.1cm} n})\preceq \overrightarrow k_{\hspace{-.1cm} 0}$, where $\preceq$ is the parity-lexicographical ordering.
\end{theorem}
Next, we equip the symbolic space with its natural  product topology. Using topological arguments, one may prove: 
\begin{theorem}
	A sequence $\bar p = \dots p_{-2} p_{-1} \sq p_0 p_1 \dots$ is admissible if and only if for every positive integer $n$ there is a $W^u$-admissible sequence $\bar q =  \dots q_{-2} q_{-1}  \sq q_0 q_1 \dots$ such that $p_{-n} \dots p_n = q_{-n} \dots q_n$.   
\end{theorem}

\begin{problem} \normalfont
	Describe the set of kneading sequences $\mathfrak{K}$. 
\end{problem}

In \cite{Is}, Ishii developed formulas that can be used to obtain a relation between parameters $a, b$, a turning point $T = (x_T, 0)$ of the Lozi map $L_{a,b}$, and its itinerary $\bar k$ (that is a kneading sequence of $L_{a,b}$). This relation is $p(a, b, \bar k) = x_T = q(a, b, \bar k)$, where $p = p(a, b, \bar k)$ is given in formula \cite[(4.2)]{I} and $q = q(a, b, \bar k)$ is given in formula \cite[(4.3)]{I}. Therefore, every kneading sequence $\bar k$ gives an equation 
\begin{equation}
	p(a, b, \bar k) = q(a, b, \bar k). 
\end{equation}
Numerical experiments show that if one has two kneading sequences, $\bar k^0$ of the rightmost turning point $T_0$ and $\bar k^{-1}$ of the leftmost turning point $T_{-1}$, and if these two turning points lie in the stable manifolds of some periodic points with small periods, then it is possible to calculate $a$ and $b$ from the corresponding two equations,
implying that these two kneading sequences govern all other kneading sequences. 

\subsubsection{Example}
\begin{figure}
	\begin{center}
		\includegraphics[height=6cm]{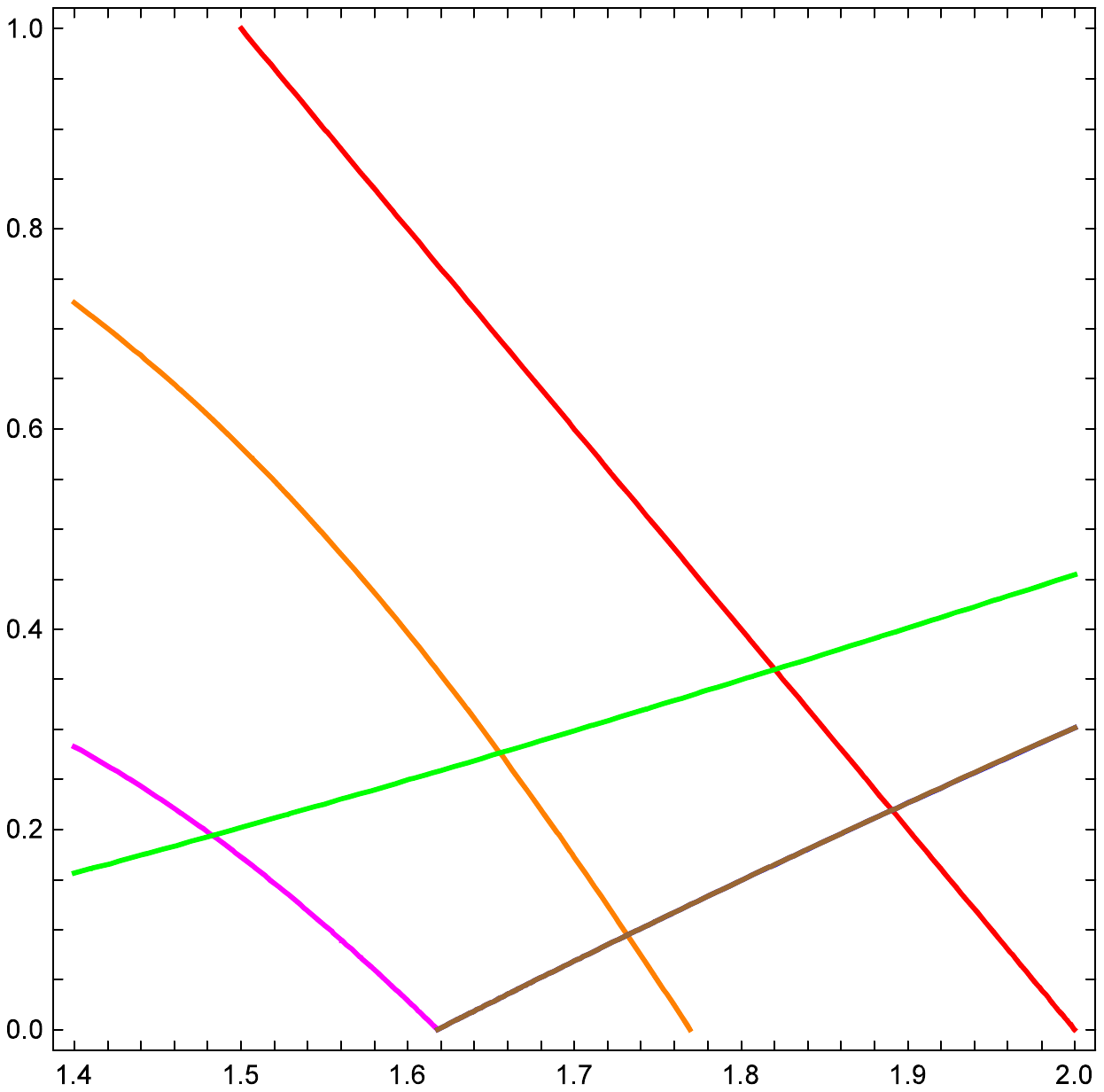}
		\caption{Graph of (2) is in orange and of (3) in green. Graph of (1) for $\bar k = \hspace{.2cm}+^{\hspace{-.5cm}\infty}\hspace{.2cm} \pm \sq (+ - -)^\infty$ is in magenta and for $\bar k = \hspace{.2cm}+^{\hspace{-.5cm}\infty}\hspace{.2cm} - \pm \sq (+ - -)^\infty$ in brown. Graph of the line $2a + b = 4$, that is a boundary line of the Misiurewicz set, is in red.}\label{fig}
	\end{center}
\end{figure}

Let $\bar k^0 = \hspace{.2cm}+^{\hspace{-.5cm}\infty}\hspace{.2cm} \pm \sq + - - +^{\hspace{-.1cm}\infty}$ and $\bar k^{-1} = \hspace{.2cm}+^{\hspace{-.5cm}\infty}\hspace{.2cm} - \pm \sq (+ -)^\infty$. The equation $p(a, b, \bar k^0) = q(a, b, \bar k^0)$ reads
\begin{equation}
	a^4 - 6a^2 -4a + 4 b^2 + a^2b + (a^3 + 2a - ab)\sqrt{4 b + a^2} = 0,
\end{equation}
and the equation $p(a, b, \bar k^{-1}) = q(a, b, \bar k^{-1})$ reads 
\begin{equation}
	\frac{4(-a^2-2b^2+2b+a \sqrt{a^2-4b})}{a-2b-\sqrt{a^2-4b}} - \left(2+a-\sqrt{a^2+4b}\right)\left(3a-\sqrt{a^2+4b}\right)=0.
\end{equation}
Using the ``NSolve'' command of Wolfram Mathematica produces a unique solution to this system of equations in the region $a\in[1,2]$, $b\in[0,1]$. This solution is approximately
$a = 1.655319602968851744592, b = 0.2765071079677260998121$, see Figure \ref{fig}.

\begin{question}\normalfont
	Is it true that any two distinct kneading sequences determine a unique pair of parameters $(a,b)$, and in that way govern all the other kneading sequences of $\mathfrak{K}$? 
\end{question}

Very recently, in \cite{BSt}, the authors developed a kneading theory for the H\'enon maps $H_{a,b}$ within a set of parameters $\mathcal{WY}$ 
for which Wang and Young proved the existence of a strange attractor. 
This set has positive measure and consist of maps which are strongly dissipative.
We refer to  \cite{BSt} and \cite{WY} for details.
\begin{problem}\normalfont
	Describe the set of kneading sequences $\mathfrak{K}$ of the H\'enon map $H_{a,b}$, with $(a,b) \in \mathcal{WY}$.  
\end{problem}	
\begin{question} \normalfont
	Is it true that any two distinct  kneading sequences of the H\'enon map $H_{a,b}$, with $(a,b) \in \mathcal{WY}$, determine a unique pair of parameters $(a,b)$, and in that way govern all the other kneading sequences in $\mathfrak{K}$ of $H_{a,b}$?  
\end{question}
\newpage


\section{H\'enon maps tangent to the identity (X. Buff)}

We propose to investigate the local dynamics of some specific H\'enon maps. 
Consider the quadratic complex H\'enon map $H_2\colon\mathbb{C}^2\to \mathbb{C}^2$ defined by 
\[H_2\left(x,y\right) = \left(y,x+y^2\right).\]
The origin is a fixed point and $H_2^{\circ 2}$ is tangent to the identity at the origin. 

\begin{figure}[htbp]
	\includegraphics[height=6cm]{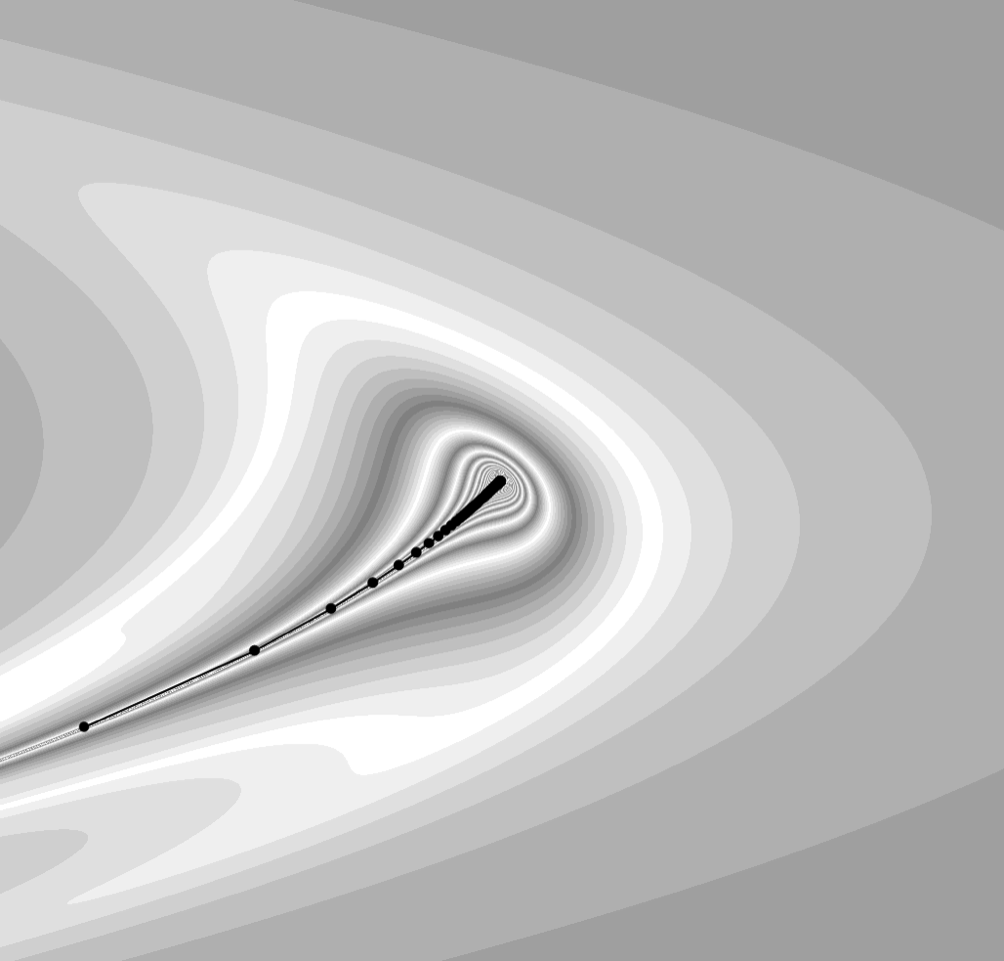}
	\includegraphics[height=6cm]{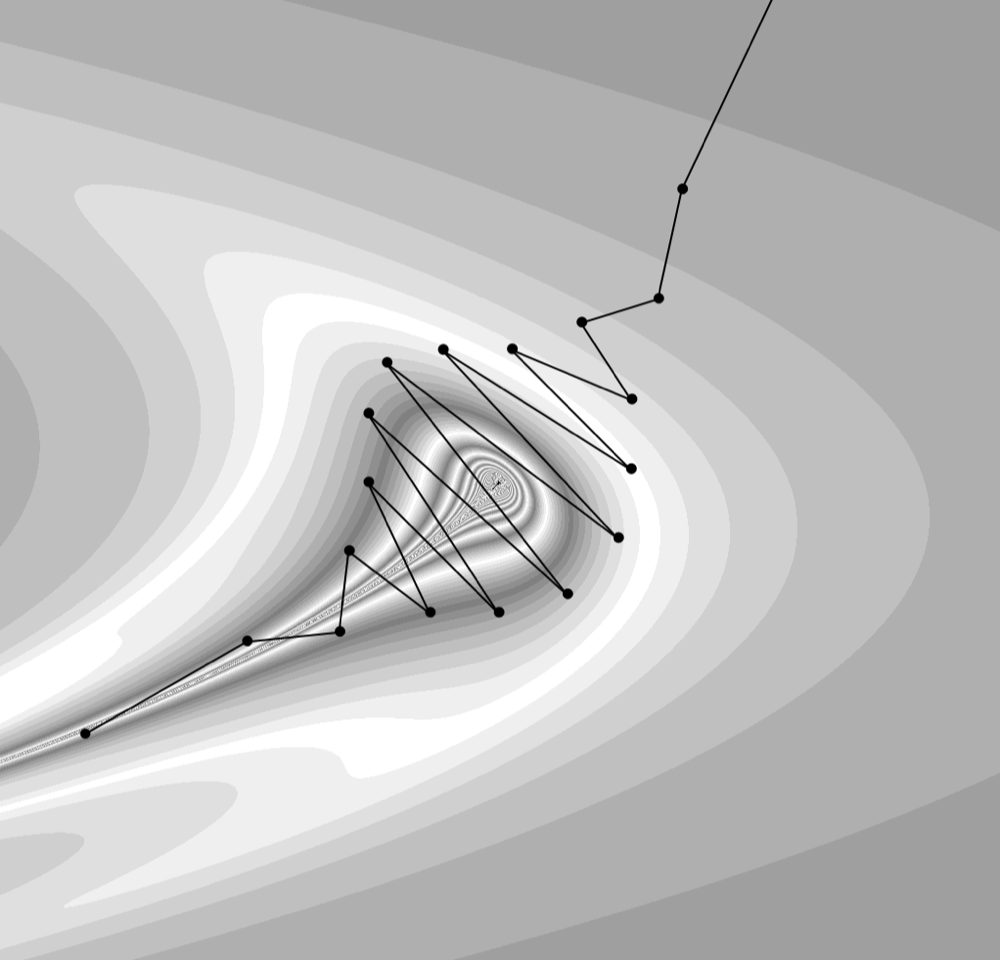}
	\caption{Left: an orbit  converging to the origin. Right: an orbit diverging to infinity.\label{fig:realH2}}
\end{figure}

Note that $H_2$ restricts to an orientation reversing diffeomorphism $H_2\colon\mathbb{R}^2\to \mathbb{R}^2$. The dynamics in $\mathbb{R}^2$ is well understood. There is an analytic map $\phi_2\colon\mathbb{R}\to \mathbb{R}^2$ such that 
\[H_2\circ \phi_2(t) = \phi_2(t+1) \ \mathrm{and} \ \phi_2(t)\sim \left(\frac{-2}{t},\frac{-2}{t}\right)\text{ as } t\to +\infty.\]
The curve $\phi_2(\mathbb{R})$ is invariant by $H_2$ and within $\phi_2(\mathbb{R})$, every orbit converges to the origin in $\mathbb{R}^2$. Outside the origin and $\phi_2(\mathbb{R})$, every orbit diverges to infinity (see Figure \ref{fig:realH2}). 

\begin{question}
	Can we describe the dynamics of $H_2$ near the origin in $\mathbb{C}^2$~?  
\end{question}

\begin{figure}[htbp]
	\includegraphics[height=6cm]{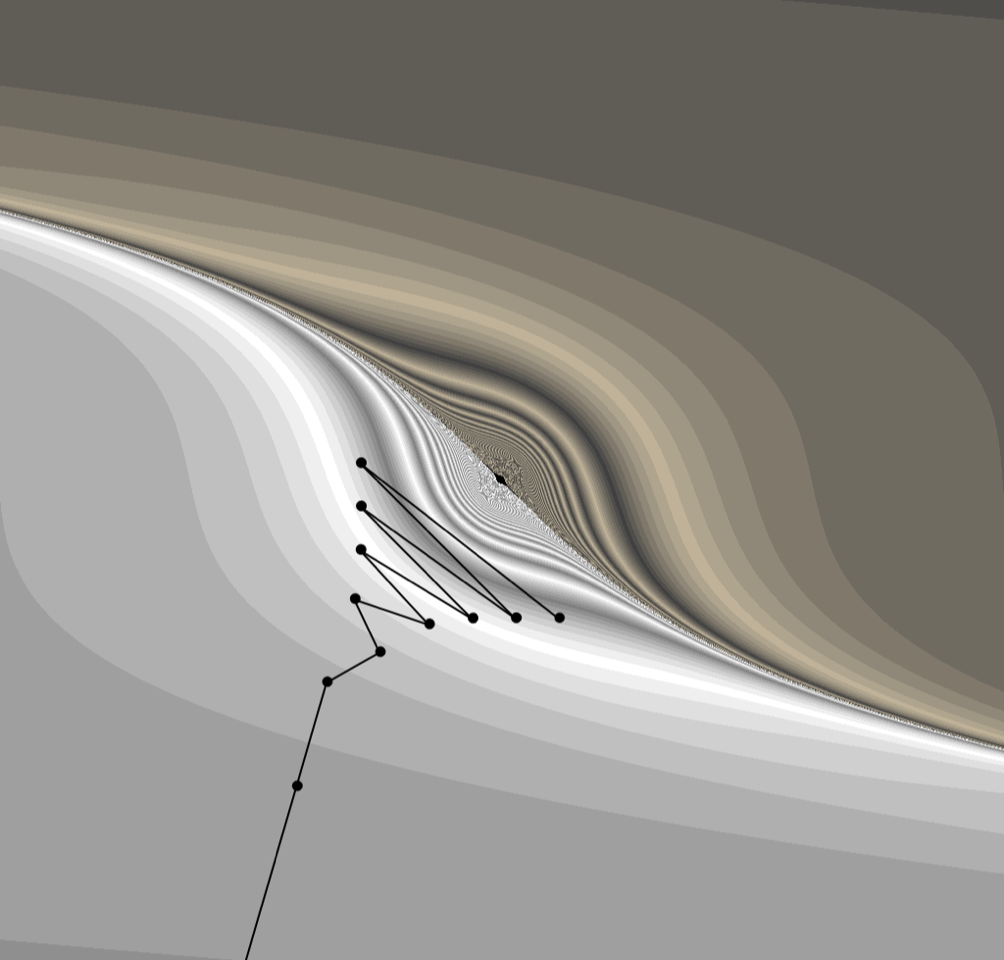}
\includegraphics[height=6cm]{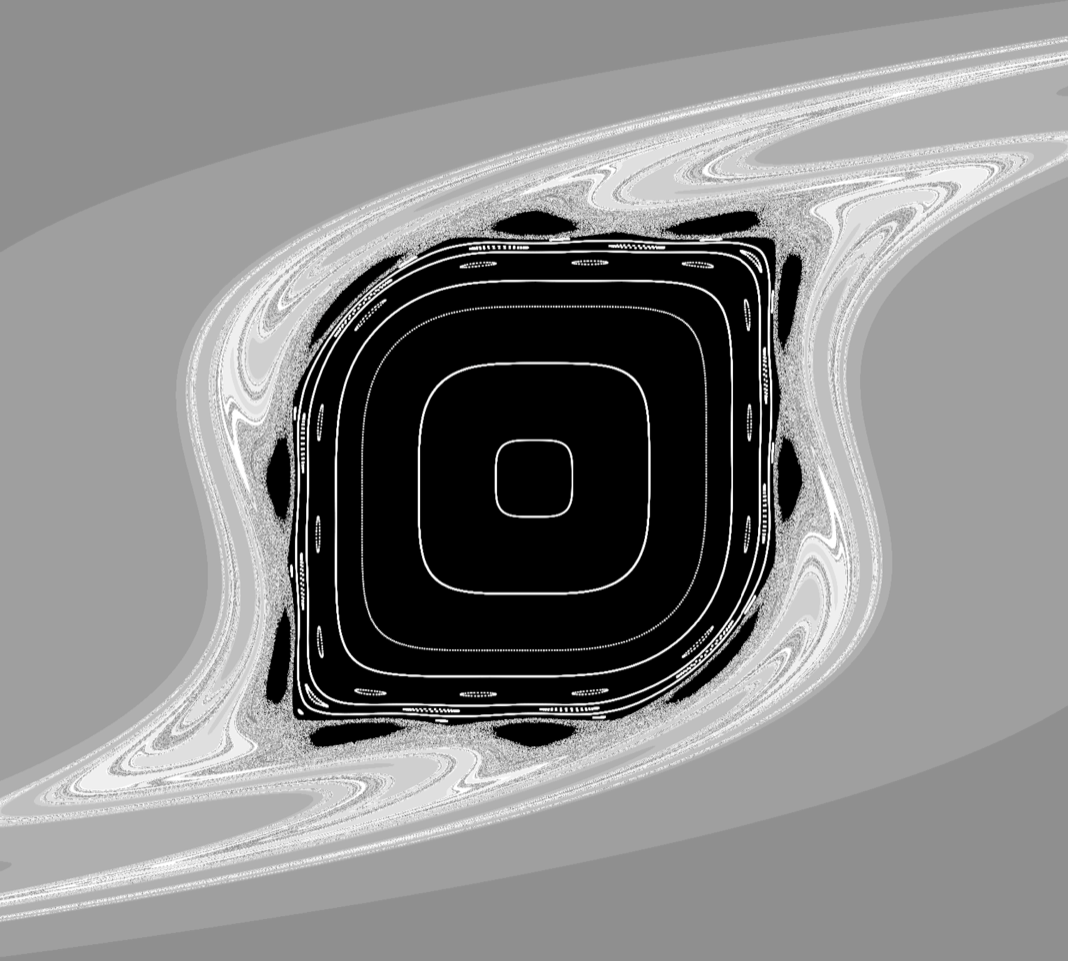}
\caption{Left: The points are colored according to whether $x+y$ tends to $+\infty$ (dark grey) or to $-\infty$ (light grey). Right: The dynamics of $H_3^{\circ 2}\colon\Pi_1\to \Pi_1$ exhibits KAM phenomena.}\label{fig:realH3}
\end{figure}

Before specifying this question, let us consider the H\'enon map $H_3\colon\mathbb{C}^2\to \mathbb{C}^2$ defined by 
\[H_3\left(x,y\right) = \left(y,x+y^3\right).\]
This H\'enon map also preserves $\mathbb{R}^2$ and the dynamics in $\mathbb{R}^2$ is also completely understood. There is an analytic map $\phi_3\colon\mathbb{R}\to \mathbb{R}^2$ such that 
\[H_3\circ \phi_3(t) = -\phi_3(t+1)\ \mathrm{and} \ \phi_3(t)\sim \left(\frac{1}{\sqrt{t}},\frac{-1}{\sqrt{t}}\right)\text{ as } t\to +\infty.\]
The curves $\phi_3(\mathbb{R})$ and $-\phi_3(\mathbb{R})$ are exchanged by $H_3$. Within those curves, every orbit converges to the origin. Outside those curves and the origin, every orbit diverges to infinity (see Figure \ref{fig:realH3} left). 

Set $\omega = e^{{\rm i}\frac{\pi}{8}}$ so that $\omega^9 = -\omega$ and consider the real planes $\Pi_1\subset \mathbb{C}^2$ and $\Pi_2\subset \mathbb{C}^2$ defined by 
\[\Pi_1 = \left\{\left(\omega x,\omega^3 y\right) ~:~x\in \mathbb{R},~y\in \mathbb{R}\right\}\ \mathrm{and} \ \Pi_2 = \left\{\left(\omega^3 x,\omega y\right) ~:~x\in \mathbb{R},~y\in \mathbb{R}\right\}.\]
Observe that $H_3$ exchanges the planes $\Pi_1$ and $\Pi_2$: 
\[H_3\left(\omega x,\omega^3 y\right) = \left(\omega^3 y, \omega (x-y^3)\right) \ \mathrm{and} \ H_3\left(\omega^3 x,\omega y\right) = \left(\omega y, \omega^3 (x+y^3)\right).\]
The dynamics of $H_3^{\circ 2}\colon\Pi_1\to \Pi_1$ is much more complex than that of $H_3\colon\mathbb{R}^2\to \mathbb{R}^2$ (see Figure \ref{fig:realH3} right).

The second iterate of $H_3$ is tangent to the identity at the origin. More precisely
\[H_3^{\circ 2}\left(x,y\right) = \left(x,y\right) + \left(y^3,x^3\right) + {\mathcal O}\bigl(\|x,y\|^4\bigr).\]
It follows that near the origin, the orbits of $H_3^{\circ 2}$ shadow the orbits of the vector field 
\[\vec{v}_3 = y^3 \partial_x + x^3\partial_y.\]
The vector field is a Hamiltonian vector field. It is tangent to the level curves of the function 
\[\Phi_3= x^4-y^4.\]
Note that
\[\Phi_3\left(\begin{array}{c}\omega x  \\ \omega^3 y\end{array}\right) = {\rm i}(x^4+y^4)\]
so that the intersection of the level curves of $\Phi_3$ with the real plane $\Pi_1$ are topological circles. Those topological circles are invariant by the flow of the vector field $\vec{v}_3$. 
It follows from the theory of Kolmogorov-Arnold-Moser that in any neighborhood of the origin, there is a set of positive Lebesgue measure of topological circles which are invariant by $H_3$ and on which $H_3$ is analytically conjugate to a rotation $\mathbb{R}/\mathbb{Z}\ni t\mapsto t+\theta\in \mathbb{R}/\mathbb{Z} $ with bounded type rotation number $\theta\in (\mathbb{R}\smallsetminus \mathbb{Q})/\mathbb{Z}$. Those invariant circles are separated by small saddle cycles and small elliptic cycles. The analytic conjugacies extend to complex neighborhoods of $\mathbb{R}/\mathbb{Z}$ in $\mathbb{C}/\mathbb{Z}$. This proves that $H_3$ has lots of Herman rings. 

Coming back to our initial problem, observe that the second iterate of $H_2$ is also tangent to the identity at the origin with
\[H_2^{\circ 2}(x,y) = (x,y) + \left(y^2,x^2\right) + {\mathcal O}\bigl(\|x,y\|^3\bigr).\]
It follows that near the origin, the orbits of $H_2^{\circ 2}$ shadow the orbits of the vector field 
\[\vec{v}_2 = y^2 \partial_x + x^2\partial_y.\]
The vector field is also a Hamiltonian vector field. 
The vector field $\vec{v}_2$ is tangent to the level curves of the function 
\[\Phi_2= x^3-y^3.\]
We can no longer apply the theory of Kolmogorov-Arnold-Moser since there is no invariant real-plane on which the level curves of $\Phi_2$ are topological circles. However, we may wonder whether 
the complex dynamics of $H_2$ exhibits KAM phenomena. 

We say that $H_2$ has small cycles if for any neighborhood $U$ of the origin $\boldsymbol{0}$ in $\mathbb{C}^2$, there exists a cycle of $H_2$ which is entirely contained in  $U\smallsetminus \{\boldsymbol{0}\}$.

\begin{question}\normalfont
	Does $H_2$ have small cycles?   
\end{question}

\begin{question}\normalfont
	Does $H_2$ have both small saddle cycles and small elliptic cycles?  
\end{question}

We say that $H_2$ has a Herman ring with rotation number $\theta\in (\mathbb{R}\smallsetminus \mathbb{Q})/\mathbb{Z}$ if there exists an annulus $V= \bigl\{z\in \mathbb{C}/\mathbb{Z}~:~{\rm Im}(z)<h\bigr\}$ with $h>0$, a holomorphic map $\phi\colon V\to \mathbb{C}^2$, and an integer $n\geq 2$ such that 
\[\forall z\in V,\quad H_2^{\circ n}\circ \phi(z) = \phi(z+\theta).\]

\begin{question}\normalfont
	Does $H_2$ have a Herman ring?  
\end{question}

If the answer is yes, we may consider the set $\Theta\subset (\mathbb{R}-\mathbb{Q})/\mathbb{Z}$ of rotation numbers $\theta$ such that $H_2$ has a Herman ring with rotation number $\theta$.  

\begin{question}\normalfont
	Does $\Theta$ have positive Lebesgue measure? More precisely, is $0$ a Lebesgue density point of $\Theta$? 
\end{question}

We believe that the answers to the previous questions are all affirmative. Regarding the following question, we do not have an opinion.  

\begin{question} \normalfont
	Assume $H_2$ has a Herman ring with bounded type rotation number $\theta$. Is it possible to find parameters $a\in \mathbb{D}\smallsetminus \{0\}$ arbitrarily close to $1$ such that the dissipative H\'enon map $H\colon\mathbb{C}^2\to \mathbb{C}^2$ defined by 
	\[H(x,y) = \left(ay,x+y^2\right).\]
	has a Herman ring with rotation number $\theta$? 
 \end{question}

\section{Quasi-hyperbolicity and uniform hyperbolicity (E. Bedford)}

\subsection{Complex H\'enon maps.}\label{sec:useful-defi}

Any H\'enon map $H_{a,P}(x,y):= ( ay + P(x), x)$ where $a\in \mathbb{C}^*$ and $P\in \mathbb{C}[x]$ is
a polynomial of degree $d\ge 2$ induces a polynomial automorphism of the affine plane $H_{a,P}\colon \mathbb{C}^2\to \mathbb{C}^2$. 

For a general polynomial automorphism $f\colon \mathbb{C}^2\to \mathbb{C}^2$, write $(f(x,y) = (P(x,y), Q(x,y))$ and
define its degree $\deg(f) := \max \{ \deg(P), \deg(Q)\}$. 
In a celebrated article~\cite{FM}, Friedland and Milnor have proved the following remarkable result (see \S \ref{sec:poly-fine} below for more details).
If the sequence of degrees $\deg(f^n)$ is unbounded, then $f$ is actually
conjugated to a composition of H\'enon maps $H_{a_1,P_1}\circ \cdots \circ H_{a_k,P_k}$, and $\deg(f) = \deg(P_1) \cdots \deg(P_k)$.

Any such composition will be  called a generalized H\'enon map.

\subsection{Quasi-expanding  H\'enon maps}\label{sec:quasi-expanding}

Suppose  $f$ is a generalized H\'enon map of degree $d\ge2$, and let ${\mathcal S}$ denote its set of (periodic) saddle points.  It is known to 
be infinite, and its distribution represents the unique measure of maximal entropy, see~\cite{BLS}.

Given a saddle point $p$, denote by $E^u_p =\{ v\in \mathbb{C}^2, \|Df^n_p v\| \to \infty\}$ the associated unstable direction,
by $W^u_\mathrm{loc}(p)$ its local unstable manifold, and by $W^u(p)$ its global unstable manifold.

We now introduce the following three conditions measuring the expansion of $f$.

\medskip

\noindent{\bf Condition 1 }  For each $p\in{\mathcal S}$, there is a metric on $E^u_p$ which is expanded by $Df_p$ with a uniform bound independent on $p$.  
More specifically, this means that there exists $\kappa > 1$ such that for each $p$, there is a metric $\| \cdot \|_j$ on $E^u_{f^j(p)}$ so that for nonzero $v\in E^u_{f^j(p)}$, we have
\[\|Df_{f^j(p)} v\|_{j+1} \geq \kappa \|v\|_j\]
for each $j$.

\medskip

If $f$ is a H\'enon map, we can consider its Green function (see~\cite{BS}) \[G^+(x,y) = \lim_{n\to\infty} \frac1{d^n} \log \max \{1, |f^n(x,y)|\}.\] We let $W_r^u(p)$ denote the connected
component of $B( p, r ) \cap W^u(p)$ containing $p$, where $W^u(p)$ is the unstable manifold at $p$ and $B( p, r )$ the euclidean ball in $\mathbb{C}^2$. \\

\noindent{\bf Condition 2 }  The unstable manifolds $W^u(p)$ satisfy the  proper, locally bounded area condition: there exist $\varepsilon > 0$ and $A < \infty$ such that for each $\delta >0$ there is
	an $\eta  > 0$ such that for each saddle point $p$ we have: $W_\varepsilon^u(p)$ is closed in $B(p, \varepsilon)$, $\mathrm{Area}(W_\varepsilon^u(p))\leq A$, and $\sup_{ W_\delta^u(p)} G^+ \geq \eta$  (see \cite[Corollary 3.5]{BS8}).

\medskip

Recall that each unstable manifold is uniformized by an entire map $\xi_p\colon \mathbb{C}\to W^u(p)\subset \mathbb{C}^2$ with $\xi_p(0)=p$.   Using the Green function, we may normalize it by putting $\hat\xi_p(\zeta):=\xi_p(\alpha \zeta)$ so that
\begin{equation}\label{eq:*}
\max_{|\zeta|\le 1} G^+(\hat\xi_p(\zeta))=1. 
\end{equation}

\medskip

The last condition we want to mention is: \\

\noindent{\bf Condition 3 }  The normalized maps $\{\hat\xi_p:p\in{\mathcal S}\}$ form a normal family of entire mappings.

A consequence of Condition 3 is that for all $x\in\overline{\mathcal S}$ we define a family of unstable manifolds by setting $W^u(x):=\xi_x(\mathbb{C})$.

\bigskip

For generalized H\'enon maps, Conditions 1, 2, and 3 are equivalent (see \cite{BS8}), and in case one/all of them hold, we say that $f$ is {\it quasi-expanding}.  A map $f$ is said to be {\it quasi-hyperbolic} if both $f$ and $f^{-1}$ are quasi-expanding.

Furthermore, a theorem  from \cite{BGS} asserts a quasi-hyperbolic map is uniformly hyperbolic if and only if there is no tangency between $W^u(x_1)$ and $W^s(x_2)$ for any $x_1,x_2$ in the closure of ${\mathcal S}$.

Recall the following standard definitions (see~\cite{BS}). Let $K^\pm$ denote the set of points with bounded forward orbits for $f^\pm$. We denote $J^\pm:= \partial K^\pm$. We also define $J:= J^+ \cap J^-$ and $J^*$ to be the closure of ${\mathcal S}$ (hence $J^* \subset J$).

\begin{question} \normalfont\label{Q17}
	If $f$ is quasi-hyperbolic, then is ${\rm int}(K^+)$ the union of a finite number of basins of sink orbits?  
\end{question}

\begin{question} \normalfont\label{Q18}
	If $f$ is quasi-hyperbolic, is $J=J^*$?  
\end{question} 

\begin{question} \normalfont\label{Q19}
	 If $f$ is quasi-hyperbolic, then is there no wandering Fatou component? 
\end{question} 

\begin{question} \normalfont\label{Q20}
	If $f$ is quasi-hyperbolic, do the unstable slices satisfy a John-type condition (as in \cite{BS7})? 
\end{question} 

\begin{question} \normalfont\label{Q21}
	If $f$ is quasi-hyperbolic and dissipative, and if $J$ is connected, do the external rays land at $J$?  Is $J$ a finite quotient of the real solenoid? 
\end{question} 

The answers to all these questions are  ``yes'' in the uniformly hyperbolic case by~\cite{BS7}. We thus ask whether these properties remain true in the quasi-hyperbolic case (in which case the dynamics is expected to be close to being hyperbolic). 
Observe that Questions~\ref{Q17} and~\ref{Q20} have negative answer for a general H\'enon map, and that Question~\ref{Q19} was disproved recently by Berger and Biebler for
some H\'enon map of degree $6$ (which are known not to be quasi-hyperbolic). 
Question~\ref{Q18} is a well-known  problem. Beside the hyperbolic case~\cite{BS}, it has been solved in a few other cases~\cite{MR3825010,MR4216565}, but the case of a generalized H\'enon map remains elusive.

\subsection{Surface automorphisms}

Let $X$ be any compact complex K\"ahler surface, and let $f\colon X\to X$ be any holomorphic
automorphism having positive topological entropy $\log \lambda$, $\lambda >1$. 
By Cantat~\cite{Cantat}, $X$ is isomorphic to either the blow-up of $\mathbb{P}^2$ at at least $10$ points, or a K$3$ surface, or an Enriques surface, or an Abelian surface.

In this context, Conditions 1 and 2 are still meaningful, but we do not have a Green function $G^+$.  However, since the dynamical degree of $f$ is $\lambda>1$, there exists an expanded positive closed $(1,1)$ current $T^+$ with $f^*T^+=\lambda T^+$ (see~\cite{Cantat}).  In this case, we can replace the normalization~\eqref{eq:*} with a condition involving the mass of a slice of the current $T^+$:
\begin{equation}\label{eq:*2}
{\rm Mass}(T^+|_{\xi_p(|\zeta|<1)}) = 1
\end{equation}
We can thus formulate a Condition $3'$, which is Condition 3  with the normalization~\eqref{eq:*2}.  

\begin{question} \normalfont
	Is Condition $3'$ equivalent to 1 and 2?  And do Questions~\ref{Q17}, \ref{Q18}, \ref{Q19} above hold for quasi-hyperbolic surface automorphisms?  
\end{question}

\subsection{Real maps}

Let us now suppose that $f$ is a real surface automorphism. In other words, we suppose $X$ to be projective and defined by real polynomial equations,  
and $f$ to be also defined over the real numbers. 
We may thus consider the restriction map $f_\mathbb{R}$ to the set of real points $X_\mathbb{R}$. 
Observe that $h_{\mathrm top} (f_\mathbb{R})\le \log \lambda$. It was proved in 
\cite{BS8} that for generalized H\'enon maps $h_{\mathrm top} (f_\mathbb{R})= \log \lambda$ implies $f$ to be
 quasi-hyperbolic.  

\begin{question} \normalfont
	 If $f$ is a real surface automorphism such that the entropy of $f_\mathbb{R}$ is the same as the entropy of $f$, does it follow that $f$ is quasi-hyperbolic?  
\end{question}

We refer to \cite{DK} for a discussion of real surface automorphisms satisfying this condition on the entropy.


\section{Parameter loci for the H\'enon family (Y. Ishii)}

\subsection{Connectedness locus}
Consider the complex H\'enon family:
\[f_{c, b} (x, y):= (x^2+c-by, x),\]
where $(c, b)\in\mathbb{C}^2$ is a parameter.\footnote{We include the case $b=0$ to simplify the presentation.} Let $J_{c, b}$ be the Julia set of $f_{c, b}$: by definition this
is the intersection between the boundaries of the sets of points having bounded forward (resp. backward) orbits.
By extension, we let $J_{c, 0}$ be the Julia set of $p(z)=z^2+c$. The \emph{connectedness locus} of $f_{c, b}$ is defined as
\[\mathcal{M}=\big\{(c, b)\in\mathbb{C}^2 : J_{c, b} \mbox{ is connected}\big\}.\]

\begin{conj}
	$\mathcal{M}$ is disconnected.  
\end{conj}

It has been shown that $\mathcal{M}\cap\mathbb{R}^2$ is disconnected~\cite{AI2}, which partially supports the conjecture above.

\subsection{Horseshoe locus}
We say that $f_{c, b}$ is a \emph{complex hyperbolic horseshoe} if $J_{c, b}$ is a hyperbolic set for $f_{c, b}$ and the restriction $f_{c, b} : J_{c, b}\to J_{c, b}$ is topologically conjugate to the full $2$-shift. The \emph{complex hyperbolic horseshoe locus} is defined as
\[\mathcal{H}_{\mathbb{C}}=\big\{(c, b)\in\mathbb{C}^2 : f_{c, b} \mbox{ is a complex hyperbolic horseshoe}\big\}.\]
One can see that $\mathcal{H}_{\mathbb{C}}$ is not simply connected since the monodromy representation:
\[\rho : \pi_1(\mathcal{H}_{\mathbb{C}})\longrightarrow \mathrm{Aut}(\{0, 1\}^{\mathbb{Z}})\]
of the fundamental group of $\mathcal{H}_{\mathbb{C}}$ (with the base-point at $(c, b)=(-4, 0)$) to the group of shift-commuting automorphisms of $\{0, 1\}^{\mathbb{Z}}$ is surjective (see, e.g.,~\cite{A,BS1}). 

\begin{question}
	Is the locus $\mathcal{H}_{\mathbb{C}}$ connected?  
\end{question}

For $(c, b)\in\mathbb{R}^2$, we can consider the restriction of $f_{c, b}$ to $\mathbb{R}^2$ and we can analogously define the \emph{real hyperbolic horseshoe locus} $\mathcal{H}_{\mathbb{R}}\subset\mathbb{R}^2$. One of the main result of~\cite{AI1} states that $\mathcal{H}_{\mathbb{R}}$ is connected and simply connected (see also~\cite{BS2}).

\subsubsection{Isentropes}
Take again $(c, b)\in\mathbb{R}^2$ and consider the restriction $f_{c, b}|_{\mathbb{R}^2} : \mathbb{R}^2\to\mathbb{R}^2$. Let $h_{\mathrm{top}}(f_{c, b}|_{\mathbb{R}^2})$ be the topological entropy of the real H\'enon map $f_{c, b}|_{\mathbb{R}^2}$. For every $0\leq \alpha \leq \log 2$, the \emph{isentrope} is defined as
\[\mathcal{E}_{\alpha}=\big\{(c, b)\in\mathbb{R}^2 : h_{\mathrm{top}}(f_{c, b}|_{\mathbb{R}^2})=\alpha\big\}.\]
In a topological term, monotonicity of the topological entropy of the real H\'enon map $f_{c, b}|_{\mathbb{R}^2}$ can be formulated as 

\begin{question}[van Strien \cite{vanstrien}]
	Is the isentrope $\mathcal{E}_{\alpha}$ connected for any $0\leq \alpha \leq \log 2$ ? 
\end{question}

Milnor and Tresser~\cite{MT} showed it is true for cubic polynomials. The main result of~\cite{AI1} implies that the locus $\mathcal{E}_{\log 2}$ is connected and simply connected (see also~\cite{BS2}).

Several articles attempt at giving lower bounds for topological entropy of real H\'enon maps, e.g.,~\cite{NBGM,NP}. Among others, the paper~\cite{NBGM} rigorously showed that $h_{\mathrm{top}}(f_{c, b}|_{\mathbb{R}^2})>0.46469$ for the classical H\'enon's parameter, and this bound is believed to be close to optimal. For upper bounds, the paper~\cite{BLS} has shown  that $h_{\mathrm{top}}(f_{c, b}|_{\mathbb{R}^2})<\log 2$ if and only if the Julia set of $f_{c, b}$ (as a complex dynamical system) is not contained in $\mathbb{R}^2$. However, there is no algorithm which provides rigorous (non-trivial) upper bounds.  

We thus propose the following problem.

\begin{question}
	Construct an algorithm to compute a rigorous upper bound for the topological entropy of a real H\'enon map $f_{c, b}|_{\mathbb{R}^2}$.  
\end{question}

Probably the only existing formula for (non-trivial) upper bound is given in~\cite{Y}. However, according to Yomdin himself, the bound in the current form is far from sharp and would not give non-trivial ones.


\section{Topology and rigidity of H\'enon maps (R. Dujardin)}

For polynomials and rational maps in dimension 1, there is a well-known list of exceptional examples 
whose Julia sets and dynamical properties are unexpectedly regular: Chebychev polynomials, monomial 
mappings and Latt\`es examples. They can characterized in many different ways, see e.g.,~\cite{MR1032883,MR2325017}.

For generalized H\'enon maps (as defined in \S\ref{sec:useful-defi}) 
it is expected that no such exceptional example exists, but not so many actual results  
in this direction are known:
\begin{itemize}
	\item Brunella proved in \cite{brunella} that a generalized H\'enon map cannot preserve an 
	algebraic foliation of $\mathbb C^2$, i.e., a singular algebraic foliation by holomorphic curves.  Here by preserving we mean that $f$ maps leaves into leaves. 
	\item   Bedford and Kim 
	proved in \cite{BK1, BK2} that
	neither $J^+$ nor $J^-$ (see \S \ref{sec:quasi-expanding} for a definition) can be a smooth $C^1$ submanifold, nor a semi-analytic set. 
\end{itemize}

\medskip

Here we propose a few rigidity questions related to these results.  

The first  question is about a quantitative reinforcement of the Bedford-Kim theorem. 
Recall from the introduction, the definition of the standard quadratic H\'enon map
$H_{a,c} (x,y) := (ay+x^2+c, a x)$, and the definition of $K^+$ and $J^+$ from \S\ref{sec:quasi-expanding}.

\medskip
 
For $(a,c)$ close to $(0, 0)$, $H_{a,c}$ is a small perturbation of the monomial map $(x,0) \mapsto (x^2, 0)$, whose 
Julia set is smooth, and in this case $J^+_{a,c}$ is a topological 3-manifold.  

\begin{question}
\normalfont	Give an asymptotic expansion of the Hausdorff dimension of $J^+_{a,c}$ as $(a,c)$ tends to $(0, 0)$. In particular is there   a uniform lower bound of $\dim(J^+_{a,c})$ of the form $\dim(J^+_{a,c})\geq 3+h(a)$ with $h(a)>0$ in the neighborhood of $c=0$?  
\end{question}

Note that $(a,c)\mapsto \dim(J^+_{a,c})$ is real analytic in the domain where $H_{a,c}$ is hyperbolic 
(this was proved for one-dimensional maps by Ruelle in~\cite{MR0684247}, and by Wolf~\cite{MR1813546} for polynomial automorphisms).  
It is not clear whether the dimension of the Julia set remains real-analytic when $H_{a,c}$ degenerates to a 
unimodal map, for instance in a full neighborhood of $(a,c) = (0,0)$.
\medskip

Can we make Brunella's theorem local? More precisely: 

\begin{question}\normalfont
	Is it possible to find a generalized complex H\'enon map $f$,
	an open set $U\subset \mathbb C^2$ intersecting $J^+$ and a holomorphic foliation of $U$ such that 
	$J^+\cap U$ is a union of leaves of this foliation?  
\end{question}

We conjecture that the answer to this question is ``no''.  The answer is presumably easier if we assume that $U\cap J^*\neq \emptyset$.  It is also possible that if $f$ is dissipative (i.e., 
$|\det(D f)|<1$), the assumption that $J^-$ is foliated is stronger than the assumption that $J^+$ is foliated (see \cite[\S 2]{bedford-dujardin}).

This would imply in particular that a generalized H\'enon map cannot preserve a (transcendental) holomorphic foliation $\mathcal F$ of $\mathbb{C}^2$. Indeed in such a case, consider the leaf $L$ through a saddle periodic point $p$:  $L$ must be mapped into itself by $f^n$, hence coincide with the stable $W^s(p)$ or unstable $W^u(p)$ manifolds (see again \S \ref{sec:quasi-expanding} for a discussion of these objects); 
changing $f$ to $f^{-1}$ if necessary, we may assume that 
$L = W^s(p)$, and since $W^s(p)$ is dense in $J^+$ it follows that $J^+$ is a 
union of leaves  of $\mathcal F$.

\begin{remark}\label{rem:fol-om}
Note that the basin of attraction of the super-attracting point at infinity $\Omega(f) := \{(x,y), |f^n(x,y)| \to \infty\}$ carries a natural 
(transcendental) holomorphic foliation which is $f$-invariant defined by the holomorphic $1$-form $\partial G^+$, see~\cite[\S 7]{HOV}.
However this foliation does not extend to  $\mathbb{C}^2$ (otherwise it would extend holomorphically to $\mathbb{P}^2$ which is absurd, see~\cite[\S 3]{MR2209087}).
\end{remark}

Related results include:
\begin{itemize}
	\item  the classification of holomorphic Anosov diffeomorphisms on surfaces 
	by Ghys \cite{ghys}, in which a basic 
	step is to prove that stable and unstable laminations are actually holomorphic foliations;
	\item the classification of birational maps preserving algebraic foliations by 
	Cantat and Favre \cite{cantat-favre};
	\item the work of Pinto, Rand and others on the smooth rigidity of hyperbolic diffeomorphisms on surfaces (see e.g.,~\cite{pinto-rand-ferreira}). 
\end{itemize}

 If the stable lamination is holomorphic, then by holonomy 
the unstable slices are holomorphically equivalent. We can now forget the foliation and ask about 
holomorphic equivalence of stable/unstable slices.  

\begin{question} \normalfont
	Under which circumstances is it possible that there exist saddle points $p$ and $q$, and  relative 
	open subsets $U\subset W^u(p)$ and $V\subset W^u(q)$ 
	such that $U\cap K^+$ is  biholomorphic to $V\cap K^+$?  
\end{question}

One obvious possibility is that $p$ and $q$ belong to the same cycle, and that 
the biholomorphism is induced by the action of $f$. We suspect that this is the only possibility. 

A variant of this problem is when $p$ and $q$ are associated to different mappings.  

\begin{question}\normalfont
	Let $f_1$, $f_2$ be two generalized H\'enon maps. Under which circumstances is it possible
	that some local unstable slice of $f_1$ (i.e., a set of the form $K^+\cap U$, where $U$ is a relatively
	open subset of an unstable manifold)
	is biholomorphically equivalent to an unstable slice of $f_2$?  
\end{question}

Since a local unstable slice of a generalized H\'enon map contains essentially complete information,
we expect that this can happen only if $f_1$ and $f_2$ 
are related by some algebraic correspondence. 
This question was raised in \cite[Remark 4.4]{DF2} for $f_2 = f_1^{-1}$, and a complete understanding would imply the main conjecture of \cite{DF2}. The analogous question of 
existence of local biholomorphisms between Julia sets 
for 1-dimensional rational maps  was addressed in \cite{symmetries, ji-xie}.  

 Since a local unstable slice of a generalized H\'enon maps contains 
essentially a complete
information about unstable multipliers, the previous question is reminiscent of the classical 
``spectral rigidity'' problem: 

\begin{question}\normalfont
	To which extent is a generalized H\'enon map determined by the list of its unstable multipliers (resp. by the list of moduli of its unstable multipliers)?  
\end{question}

We refer to~\cite{ji-xie_multipliers} for a proof that the list of moduli of all multipliers
determine a finite set of conjugacy classes of rational map of the Riemann sphere. 


\section{Statistical properties of complex H\'enon maps (F. Bianchi and T.-C. Dinh)}

We denote in this section by $f$ a complex H\'enon map and by $\mu$ its unique measure of maximal entropy \cite{BS, BLS, Sibony99}. We are interested in the statistical properties
of $\mu$ and of other
 natural invariant measures associated to such systems.

\subsection{Thermodynamics for H\'enon maps}

Consider a continuous function $\phi\colon \mathbb C^2 \to \mathbb R$, that will be called 
a \emph{weight}. 
Following \cite{Ruelle78}
one can define
the \emph{pressure} $P(\phi)$ as
\[P(\phi) := \sup \big(h_\nu + \langle \nu,\phi \rangle\big),\]
where the supremum is taken over all invariant probability measures
$\nu$
for $f$ and $h_\nu$ denotes the 
measure-theoretic entropy of $\nu$. 
A measure $\nu_0$ maximising the above supremum is called an \emph{equilibrium state} associated to $\phi$ and is necessarily ergodic when it is unique. 
The equilibrium state associated to 
$\phi\equiv 0$ 
is the measure of maximal entropy $\mu$.
We refer to \cite{PU10} for an account on the properties of equilibrium states in one-dimensional complex dynamics
and to \cite{Baladi00,CRH23,CPZ20} and references therein
for the case of real H\'enon maps and
diffeomorphisms of compact manifolds satisfying some hyperbolicity assumptions.

\begin{problem}\label{pb_equilibrium}
	Prove the existence and the uniqueness of
	the equilibrium state $\mu_\phi$ 
	associated to any 
	sufficiently regular weight
	$\phi$
	(for instance, every H\"older continuous $\phi$, 
	and possibly with some bound on $\max \phi - \min \phi$).
\end{problem}

Recall
that saddle points are equidistributed with respect to 
the measure of maximal entropy
\cite{BLS2}. Namely, we have 
\begin{equation}\label{eq:convergence-saddle}
	\frac{1}{d^n} \sum_{x \in SP_n} \delta_{x}\to \mu,
\end{equation}
where $d$ is the algebraic degree of $f$ (or, equivalently, $\log d$ is the topological entropy of $f$, and the 
measure-theoretic entropy of $\mu$) and $SP_n$ is the set
of the saddle $n$-periodic points of $f$.

\begin{question}
Suppose $\phi$ is sufficiently regular, so that
the equilibrium state $\mu_\phi$ exists and is unique. 	
	Is it true that
	\begin{equation}\label{eq:convergence-saddle-weights}
		\frac{1}{e^{n P(\phi)}} \sum_{x \in SP_n} e^{S_n (\phi)} \delta_x \to \mu_\phi \quad ?
	\end{equation}
\end{question}

A version of the previous question has been established 
in \cite{BD23-p1} in the 
(expanding) setting of endomorphisms of $\mathbb P^k_{\mathbb C}$, and in particular for polynomials maps
on $\mathbb C$.

\medskip

Of a somehow different flavour, we recall
that an explicit speed of convergence in \eqref{eq:convergence-saddle}
is unknown. We believe  it is a very natural and challenging question to
quantify such convergence when testing against sufficiently regular functions.

\begin{question}
	Is the convergence \eqref{eq:convergence-saddle} exponentially fast against H\"older continuous observables? Is that
	also the case for \eqref{eq:convergence-saddle-weights}?
\end{question}

\subsection{Statistical properties of equilibrium states and spectral gap for the transfer operators}

Suppose the existence and the uniqueness of an equilibrium state $\mu_\phi$ have been established. 
The (deterministic)
 problem of describing all orbits in the support of $\mu_\phi$ is essentially impossible as this support should be contained
 in the set of points with chaotic behaviour in both forward and backward time. It is natural to adopt
a probabilistic (or statistical) approach to this problem,  to consider an \emph{observable} $g\colon \mathbb C^2 \to \mathbb R$, and
 to view the sequence $\{g\circ f^j\}_{j \in \mathbb N}$ as a sequence of
 random variables 
 on the probability space $(\mathbb C^2, \mu_\phi)$.
Since  $\mu_\phi$ is
 invariant,
 these random 
 variables have the same distribution. They are however not independent, since
they arise
from a deterministic setting. 
The first goal is thus to show that the correlations $\langle \mu_\phi, g\circ f^{j_1} \cdot g\circ f^{j_2} \rangle - \langle \mu_\phi, g \rangle^2 $ 
go to zero in a quantifiable way, as $|j_2-j_1|\to \infty$,
see for instance \cite[Problem 2]{Viana97}.
When this happens and the convergence is fast enough, the sequence 
$\{g\circ f^j\}_{j \in \mathbb N}$ is then expected to satisfy 
a list of properties which 
are typical of  independent identically distributed (i.i.d.) random variables.

\medskip

As a first step, since $\mu_\phi$ is ergodic, 
Birkhoff 
theorem asserts that
\begin{equation}\label{eq:birkhoff}
\frac{1}{n} S_n (g)  (x):= \frac{1}{n} (g (x) + g\circ f (x) + \dots + g\circ f^{n-1} (x)) \to \langle \mu_\phi, g \rangle := \int_{\mathbb C^2} g\, {\rm d}\mu_\phi
\end{equation}
for $\mu_\phi$-almost every $x$ and every $g\in L^1 (\mu_\phi)$. This can be seen as a version of the law of large numbers in this setting. The next step is to show the Central Limit Theorem (CLT) for sufficiently regular observables. 
As in the case of i.i.d. random variables, this CLT gives
the rate of the above convergence \eqref{eq:birkhoff}.
\begin{problem} Show that $\mu_\phi$ satisfies the CLT for H\"older continuous observables. Namely, prove that, for any H\"older continuous observable $g$, there exists $\sigma \geq 0$ such that 	for any interval $I \subset \mathbb R$, we have
	\[
	\lim_{n\to \infty} \mu_\phi \left(\left\{
	\frac{S_n (g) - n\langle \mu_\phi,g\rangle }{\sqrt{n} } \in I
	\right\}\right)
	=
	\begin{cases}
		1 \mbox{ when } I \mbox{ is of the form } I=(-\delta,\delta)  & \mbox{ if } \sigma^2=0,\\
		\frac{1}{\sqrt{2\pi\sigma^2}}\displaystyle\int_I e^{-t^2 / (2\sigma^2)}\, {\rm d}t &  \mbox{ if }
		\sigma^2 >0.
	\end{cases} \]
\end{problem}

\medskip
In the case of the maximal entropy measure $\mu$, the CLT was established in \cite{BianchiDinh1}. A natural question is also to characterize the observables for which $\sigma =0$. 

\medskip

Sequences of (almost)
independent random variables are also expected to satisfy \emph{large deviations} properties.  Recall that a coboundary $g$ is an observable of the form $\varphi \circ f - \varphi$.  

	\begin{problem}\label{question_LDP}
Show that $\mu_\phi$ satisfies the Large Deviation Principle (LDP)	 for H\"older continuous observables. 
	Namely, prove that,  for any  H\"older continuous observable   $g$ with $\langle \mu_\phi, g\rangle=0$ and which is not a coboundary, there exists a
	non-negative, strictly convex function $c$ which is defined on a neighborhood of $0 \in \mathbb R$,
	vanishes
	only at $0$, and such that, for all
	$\epsilon > 0$ sufficiently small,
	\[
	\lim_{n\to \infty}
	\frac{1}{n}
	\log \mu_\phi
	\left(\left\{
	x \in X\colon
	\frac{S_n(g)(x) }{n}
	>\epsilon
	\right\}\right)
	= -c(\epsilon).\]
	\end{problem}
 Note that this question is still open even for the measure of maximal entropy.

\medskip

A possible unified 
approach to the above statistical properties would be to find
Banach spaces
(containing all H\"older continuous functions)
where a suitable \emph{Ruelle-Perron-Frobenius (transfer) operator} 
associated to $f$ would turn out to be a strict contraction
on the complement of an invariant line, see for instance \cite{Baladi00, Gou15,Ruelle78}.
In the case of 
endomorphisms of $\mathbb P^k_{\mathbb C}$ in any dimension (in particular
for any $1$-dimensional complex polynomial), this is the main result of \cite{BD22} (see also \cite{MS00,Ruelle92}).
\begin{question}\label{37}
	Do there exist norms
	for functions on the Julia set which are bounded on H\"older continuous functions, contract
	(on the complement of an 
	invariant line) under the action of $f_*$ (or, more generally, of $f_* (e^{\phi-P(\phi)} \cdot)$),
	and such that the contraction is stable by perturbation of $\phi$?
\end{question}

In the case of hyperbolic maps, 
such a good Banach space
has been introduced by Blank-Keller-Liverani
\cite{BKL,Liverani95}.
The norm is 
obtained by 
combining
a regularity condition on the unstable manifolds together with a dual condition on the stable manifolds. Note that this was the starting point of a long story (see, e.g., \cite{GL06, BT07}).
As H\'enon maps are only non-uniformly hyperbolic
(so that stable and unstable manifolds do not behave nicely in general)
the above result does not apply here.

\medskip

A positive answer to Question~\ref{37}  would also give a unified proof for many  statistical properties of independent interest (including the Large Deviations as in Problem~\ref{question_LDP}), without the need of an \emph{ad hoc} proof for each of them. For instance, the Local Central Limit Theorem (LCLT) and the Almost Sure Invariance Principle (ASIP) are both satisfied by sequences of i.i.d., and provide stronger results than the CLT, see, e.g., \cite{PS75,Gou15} for definitions and criteria.

\begin{problem}\label{Problem_LCLT} Let $\mu_\phi$ as in Problem~\ref{pb_equilibrium}. 
	Show that all H\"older continuous observables which are not coboundaries satisfy
	the LCLT and the ASIP
	with respect to $\mu_\phi$.
\end{problem}
Recall that the ASIP implies the 
Almost Sure Central Limit Theorem 
and the Law of the Iterated Logarithm.

\subsection{Higher dimension and other generalizations}

Until now,
we restricted
our attention
 to H\'enon maps, i.e., polynomial automorphisms of $\mathbb C^2$.
On the other hand,
Problems~\ref{question_LDP}  and \ref{Problem_LCLT}, and  Question~\ref{37} make perfect sense for the equilibrium measures  of general rational maps once this measure has been successfully defined. 
We review below some partial results that have been obtained in more general (invertible) settings than Hénon maps. 

\medskip

A polynomial automorphism of $\mathbb C^k$ is said to be regular 
if the indeterminacy sets of the extensions to $\mathbb P^k_{\mathbb C}$ of $f$ and $f^{-1}$
are non-empty and disjoint (observe that every H\'enon map in dimension 2 satisfies this assumption, as these two sets are two distinct
points).
The construction of the measure of maximal entropy is given in \cite{Sibony99}, and 
the equidistribution of saddle points with respect to this measure is proved in \cite{Dinh_Sibony_equi} (see
\cite{BianchiDinh1,Dinh_decay}
for the exponential mixing and the CLT in this case).
More generally,
 one can also consider birational meromorphic
 maps of $\mathbb P^k_{\mathbb C}$, see
  \cite{BD05,DF01, Dethelin_Vigny, Dujardin_duke}
  for the construction of the measure of maximal entropy 
 and its properties.

\medskip

 Given integers 
$1\leq p<k$ and 
open bounded convex domains
 $M \Subset \mathbb C^p$  and
$N\Subset \mathbb C^{k-p}$,
a \emph{horizontal-like map} is 
a proper holomorphic map from a vertical subset to a horizontal subset of $M\times N$
which geometrically 
expands in $p$ directions and contracts in $k-p$ directions, see 
\cite{Dinh_Sibony_horizontal} for the precise definition.
 In this setting, the unique measure of maximal entropy has been constructed and studied in 
 \cite{DNS08,Dinh_Sibony_horizontal,MR2045508}.
In the invertible case,
  the CLT for this measure can be deduced from \cite{BianchiDinh1}.

\medskip
One can also consider automorphisms of compact K\"ahler manifolds, see for instance 
\cite{Cantat, DS05, Dujardin_duke}
 for the construction of the measure of maximal entropy and its properties.
This setting
shares a number of features with that of H\'enon maps (in dimension 2) and regular automorphisms (in any dimensions).
On the other hand, the compactness of the manifold makes 
it more difficult to apply pluripotential techniques as in the case of H\'enon maps. For instance, the proof of the CLT for the measure of maximal entropy,
given in 
\cite{BD23},
 requires the use of the theory of \emph{superpotentials} on such manifolds
\cite{Dinh_Sibony_Kahler}.

\section{Towards higher dimensions and complex differential geometry (C. Favre)}

\subsection{H\'enon maps and the group of polynomial automorphisms of $\mathbb{C}^2$.}\label{sec:poly-fine}

Let  $\mathrm{Aut}[\mathbb{C}^2]$ be the group of polynomial automorphisms of $\mathbb{C}^2$. 
Recall the definition of degrees of a polynomial automorphism of the affine plane from \S\ref{sec:useful-defi}. 
Jung~\cite{jung} proved that the group $\mathrm{Aut}[\mathbb{A}^2_\mathbb{C}]$
is generated by affine transformations and triangular maps of the form $E_P(x,y):= (x, y + P(y))$. 
And the more precise version of Friedland-Milnor's main theorem (\cite{FM}) states that  
either $f \in \mathrm{Aut}[\mathbb{A}^2_\mathbb{C}]$ is conjugated to a generalized H\'enon maps
$H_{a_1,P_1}\circ \cdots \circ H_{a_k,P_k}$ and $\deg(f^n) \asymp (d_1\cdots d_k)^n$ for all $n$; 
or $\deg(f^n)$ remains bounded and $f$ is conjugated to an affine map or to a triangular map. 

Let us state the following general problem in vague terms.  
\begin{problem} \normalfont
	Describe the growth type of the sequence $\{\deg(f^n)\}$ for any polynomial automorphism $f$ of  $\mathbb{C}^d$, $d\ge 3$.   
\end{problem} 
Very few results are known. Recall that Russakovski and Shiffman \cite{rus-shi} observed that 
\[
\deg(f^{n+m}) \le \deg(f^n) \deg(f^m)
\]
for all $n,m\ge0$ so that the following limit $\lambda(f) := \lim_n \deg(f^n)^{1/n}$ exists. We refer to it as the dynamical degree of $f$.

The case of cubic automorphisms on $\mathbb{C}^3$, and the case of automorphisms obtained as a composition of an affine transformation and a triangular one were considered by 
Blanc and Van Santen~\cite{blanc-V1,blanc-V2}. Their computations lead them to  formulate the following intriguing problem. 
A weak Perron number is an algebraic integer $\lambda \ge1$ such that all its Galois conjugates satisfy $|\mu|\le \lambda$.  
\begin{question} \normalfont
	Is the dynamical degree of any polynomial automorphism of $\mathbb{C}^d$ equal to a weak Perron number of degree 
	$\le d-1$?  
\end{question} 
It has been proven in~\cite{DF20}, that $\lambda(f)$ is an algebraic number of degree $\le 6$ for any polynomial automorphism $f$ of  $\mathbb{C}^3$. 

The case $\lambda(f)=1$ is also particularly interesting.  
\begin{question} \normalfont
	Suppose  $f$ is a  polynomial automorphism of  $\mathbb{C}^d$ satisfying $\lambda(f)=1$. 
	Is it true that $\deg(f^n) \asymp n^k$ for some $k \in \mathbb{N}$? Moreover, 
	if $k\ge 1$, does $f$ preserve a rational fibration?  
\end{question} 
Urech proved that $\deg(f^n)$ tends to infinity whenever it is unbounded, see~\cite{urech}. 
His result was made stronger by Cantat and Xie in~\cite{CX}: there exists a universal function $\sigma\colon \mathbb{N}\to\mathbb{N}$
such that $\limsup \sigma = \infty$ and $\deg(f^n) \ge \sigma(n)$. 
They raised the following weaker form of the previous problem.  
\begin{question} \normalfont
	Suppose that $\lambda(f)=1$, and $\deg(f^n)$ is unbounded
	for some polynomial automorphism $f$ of  $\mathbb{C}^d$. 
	Does there exist $C>0$ such that $\deg(f^n) \ge Cn$?     
\end{question}
Let $\mathrm{Tame}(3)$ be the subgroup of polynomial automorphisms  of $\mathbb{C}^3$ which is generated by affine and triangular transformations. 
A theorem of  Shestakov and Umirbaev~\cite{SU} states that $\mathrm{Tame}(3)$  is a strict subgroup 
of the full group of polynomial automorphisms  of $\mathbb{C}^3$ (as opposed to the $2$-dimensional situation). 
A decisive progress on the structure of $\mathrm{Tame}(3)$ was recently 
made by Lamy and Przytycky~\cite{LP}, who constructed a $\mathrm{CAT}(0)$-complex $\mathcal{C}$ over which $\mathrm{Tame}(3)$
acts by isometries.  

\begin{question} \normalfont
	Is it possible to characterize those  $f\in \mathrm{Tame}(3)$ for which $\lambda(f)=1$ in terms of their action on $\mathcal{C}$? 
\end{question} 
 
\subsection{H\'enon maps and compact complex manifolds.}
Consider any generalized H\'enon map
$f= H_{a_1,P_1}\circ \cdots \circ H_{a_k,P_k}$
as in the previous section. 
Recall that $f$ extends to the projective plane $\mathbb{P}^2_\mathbb{C}$ as a birational map
contracting the line at infinity to the super-attracting fixed point $p= [1:0:0]$.
The topology of the basin of attraction of this point $\Omega(f) := \{ q \in \mathbb{C}^2, \, f^n(q) \to p\}$
has been explored by Hubbard and Oberste-Vorth~\cite{HOV}.
They also observed that the map $f$ acts properly discontinuously on $\Omega(f)$ so that 
the space of orbits $S(f):= \Omega(f)/\langle f \rangle$ is naturally a complex surface. 
One can then construct a compact complex surface $\check{S}(f)$
having an isolated  normal singularity at a point $0\in \check{S}(f)$ such that 
$\check{S}(f)\setminus\{0\}$ is biholomorphic to $S(f)$.
The minimal resolution of $\check{S}(f)$ is a compact complex surface $\bar{S}(f)$
that is non-K\"ahler, contains no smooth rational curve of self-intersection $-1$, 
and satisfies $b_1(\bar{S}(f))= 1$. In Kodaira's classification of surfaces~\cite{BHPVdV}, it belongs to the class 
VII$_0$ which is arguably the most mysterious class of compact complex surfaces. 
Dloussky and Oeljeklaus~\cite{DO} have investigated when these surfaces carry global holomorphic  vector fields.

\begin{question} \normalfont
	Let $f$  be any generalized H\'enon map. 
	Describe the set of all generalized H\'enon maps $g$ such that $\bar{S}(f)$ is biholomorphic to $\bar{S}(g)$.
\end{question}  
It follows from~\cite[Proposition~2.1]{favre} that under the preceding assumptions, $f$ and $g$ have the same degree
and the same jacobian.  Some partial results have been obtained by R. Pal for maps of the form $H_{a,P}$~\cite{ratna}
generalizing former works on quadratic H\'enon maps by Bonnot-Radu-Tanase~\cite{MR3722291}.

\medskip

Surfaces $\bar{S}(f)$ carry only finitely many rational curves that are all contracted to the singular point
$0\in \check{S}(f)$. One can also prove that it carries a unique holomorphic foliation which is induced by the Levi flats
of the Green function $G^+$ on $\Omega(f)$, see Remark~\ref{rem:fol-om}.

An interesting feature of the complex surface $\bar{S}(f)$ is that it admits a family of charts $(U_i,\phi_i)$
where $U_i$ is an open cover of  $\bar{S}(f)$, and $\phi_i \colon U_i \to \mathbb{C}^2$ is an open immersion such that
$\phi_{ij}$ is the restriction to an open domain of a birational self-map of $\mathbb{P}^2_\mathbb{C}$. 

A complex manifold which admits a holomorphic atlas whose transition maps are restriction of birational maps
of $\mathbb{P}^d_\mathbb{C}$ is said to carry a birational structure. 

The following problem is extracted from~\cite{dloussky}. 
\begin{question} \normalfont
	Does any non-K\"ahler compact complex surface admit a birational structure?
\end{question}  
This question is extremely challenging, and reduces to the case of VII$_0$ surfaces. 
One can ask whether any deformation $S$ of a surface $\bar{S}(f)$ 
associated to a polynomial automorphism $f$ as above admits a birational structure. 
This is true when the surface satisfies $b_2(S) \le 3$, see~\cite{dloussky}. 

\medskip

Analogs of the construction of $\bar{S}(f)$ have been 
explored by Oeljeklaus and Renaud in~\cite{OeR} for some quadratic polynomial automorphisms of $\mathbb{C}^3$, and further expanded
by Ruggiero~\cite[Chapter~4]{thesisMR}. 
A polynomial automorphism $f\in \mathrm{Aut}[\mathbb{C}^3]$ is said to be regular if the indeterminacy locus $I(f)$ of its extension to $\mathbb{P}^3_\mathbb{C}$
is disjoint from $I(f^{-1})$. This notion was introduced by Sibony in~\cite{panorama}.   
Let $\Omega(f)$ be the basin of attraction of $I(f^{-1})$: this is an open $f$-invariant set over which $f$
acts properly discontinuously. As above, denote by $S(f)$ the quotient space $\Omega(f)/\langle f \rangle$.  

\begin{problem} \normalfont
	Let  $f\in \mathrm{Aut}[\mathbb{C}^d]$ be any regular polynomial automorphism. 
	\begin{enumerate}
		\item
		Prove that one can find a compact complex manifold $\bar{S}(f)$ and an open immersion
		$S(f) \subset \bar{S}(f)$ such that the complement $\bar{S}(f)\setminus S(f)$ is a divisor. 
		\item
		Prove that $\bar{S}(f)$ is unique up to bimeromorphism. 
		\item 
		Describe complex objects on $\bar{S}(f)$ (analytic subvarieties, vector fields, holomorphic foliations, positive closed currents,...). 
		Compute its deformation space.  
\end{enumerate}
\end{problem}

It is unclear how to extend this construction to a larger class of polynomial automorphisms of $\mathbb{C}^3$. 
However when $\lambda(f)^2>\lambda(f^{-1})$ an invariant valuation on the ring of polynomial functions in three variables\footnote{One way to interpret geometrically such an object
is to say that it picks an irreducible subvariety in \emph{any} algebraic compactification of  $\mathbb{C}^3$ in a compatible way.}
 is known to exist by~\cite{DF20}, which suggests the following question.  
\begin{question} \normalfont
	Let  $f\in \mathrm{Aut}[\mathbb{C}^3]$ be any polynomial automorphism such that $\lambda(f)^2>\lambda(f^{-1})$. 
	Prove the existence of a projective compactification $X$ of $\mathbb{C}^3$ such that 
	the induced birational map $f\colon X \dashrightarrow X$ admits a super-attracting
	fixed point $p$ on the divisor at infinity.  
\end{question} 

Once such a compactification has been found, one can consider the basin of attraction $\Omega$ of the point $p$
and try to construct a compactification of the space of $f$-orbits in $\Omega$ as above.

\section{H\'enon maps over number fields (P. Ingram)}
Consider a sufficiently large field $k$, for example a number field. In general, one should expect to be able to construct  H\'{e}non maps of the form
\[f(x, y)=(y, F(y)-\delta x)\]
with cycles of length $\deg(F)+2$. Write $F(y)=a_0+\cdots + a_dy^d$, and let $y_0, ..., y_{d+1}$ be variables ranging over $k$. Then, insisting that $f$ sends
\begin{equation}\label{eq:cycle}(y_0, y_1)\to (y_1, y_2)\to \cdots \to (y_{d}, y_{d+1})\to (y_{d+1}, y_0)\to (y_0, y_1)\end{equation}
is the same as insisting that
\[a_0+a_1y_n+\cdots + a_dy_n^d-\delta y_{n-1}=y_{n+1},\]
for all $n\;(\mathrm{mod}\ d+2)$. One checks that the determinant of the associated Vandermonde-like matrix
\[\begin{pmatrix}
	1 & y_0 & \cdots & y_0^d & -y_{d+1}  \\
	1 & y_1 & \cdots & y_1^d & -y_{0}  \\
	\vdots & \vdots & \ddots & \vdots & \vdots \\  
	1 & y_{d+1} & \cdots & y_{d+1}^d & -y_{d} 
\end{pmatrix}
\]
is not identically zero (e.g., substituting $y_{d+1}=y_d$ into this determinant gives $\pm (y_d-y_{d-1})\prod_{0\leq i<j\leq d}(y_j-y_i)$, which is itself not identically zero), and so this matrix is invertible over some affine-open subset of $k^{d+1}$. Here, one can find coefficients $a_0, ..., a_d$, and $	\delta$ of $f$ which enact~\eqref{eq:cycle}.

Recently, Hyde and Doyle~\cite{hd} exhibited single-variable polynomials over number fields with more preperiodic points than this sort of naive interpolation construction gives. One might ask if similar phenomena could be exploited for generalized H\'{e}non maps. 

\begin{question}\normalfont
	Over a number field $K$, is it possible to construct infinite families of generalized H\'{e}non  maps of algebraic degree $d$, and $K$-rational cycles of length at least $d+3$?	Can one construct maps with $N_d$ periodic points, where $N_d- d\to \infty$, or even $N_d/d\to \infty$ as $d\to\infty$? 
\end{question}

Next, it is natural to ask about bounds in the other direction. In analogy to the Uniform Boundedness Conjecture of Morton and Silverman~\cite{Morton_Silverman}, it is natural to posit the following, in which $\operatorname{Per}(f)$ is the set of periodic points of $f$ over the algebraic closure of $K$. 
\begin{conj}\label{conj:ubc} \normalfont
	Let $K$ be a number field, let $B\geq 1$, and let $d\geq 2$. Then as $f$ varies over generalized H\'{e}non maps of degree $d$ over $K$, the quantity
	\[\#\{P\in \operatorname{Per}(f):[K(P):K]\leq B\}\]
	is bounded just in terms of $[K:\mathbb{Q}]$, $d$, and $B$.	 
\end{conj}
We have already seen why this bound must depend on $d$, and adjoining periodic points of $f$ to the base field shows why the bound must depend on $[K:\mathbb{Q}]$. As a starting point for further exploration, we mention two more readily falsifiable conjectures.

\begin{conj}[see~\cite{I}] \normalfont
	Over $\mathbb{Q}$, $(x, y)\mapsto (y, y^2+c+x)$ has no point of period $N$, other than $N\in \{1, 2, 3, 4,  6, 8\}$. 
\end{conj}

\begin{conj} \normalfont
	For all but finitely many $\delta\in \mathbb{Q}$, the $\mathbb{Q}$-rational periodic points of any $f(x, y)=(y, y^2+c-\delta x)$ with $c\in \mathbb{Q}$ have period dividing 2.   
\end{conj}
It should be noted that some infinite families of examples show that both of these conjectures, if true, would be sharp.

It seems reasonable to posit something even stronger than Conjecture~\ref{conj:ubc}. 
Write $\log^+ x=\log\max\{1, x\}$ for $x\in \mathbb{R}^+$, and for an absolute value $|\cdot|_v$, set 
\[\|x_1, ..., x_m\|_v=\max\{|x_1|_v, ..., |x_m|_v\}.\] Recall that a number field $K$ comes equipped with a standard set $M_K$ of absolute values, and we define the \emph{naive Weil height} of $P\in \mathbb{A}^N(K)$ to be
\[h(P)=\sum_{v\in M_K}\frac{[K_v:\mathbb{Q}_v]}{[K:\mathbb{Q}]}\log^+\|P\|_v.\]
Kawaguchi~\cite{K} constructed a \emph{canonical height} $\hat{h}_f$ associated to a generalized H\'{e}non map, which differs by a controllable amount from the naive height, and interacts favorably with the dynamics of $f$, satisfying for example
\[\hat{h}_f\circ f + \hat{h}_f\circ f^{-1}=\left(d+\frac{1}{d}\right)\hat{h}_f,\]
and $\hat{h}_f(P)=0$ if and only if $P$ is periodic. In light of partial results in this direction~\cite{I}, it seems reasonable to conjecture the following strengthening of Conjecture~\ref{conj:ubc}. 

\begin{conj}\normalfont
	Let $K$ be a number field, let $B\geq 1$, and let $d\geq 2$. Then, there exist an $\epsilon>0$ and constant $C$ (depending on these data) such that, as $f$ varies over generalized H\'{e}non maps of degree $d$ over $K$, the quantity
	\[\#\{P\in (\overline{K})^2 :[K(P):K]\leq B\text{ and }\hat{h}_f(P)<\epsilon h(f)+C\}\]
	is bounded uniformly, where $h(f)$ is the height of the tuple of coefficients of $f$.  
\end{conj}

Finally, let $f$ be a generalized H\'{e}non map defined over a number field $K$ with good reduction away from some finite set $S$ of primes (so, the coefficients are $S$-integers, and $a_d$ and $\delta$ are $S$-units), and let $P_0\in (\mathcal{O}_{K, S})^2$ be some non-periodic point. Set $P_{n+1}=f(P_n)$, for $n\geq 0$, and 
\[\mathfrak{a}_n = \gcd(x(P_n)-x(P_0), y(P_n)-y(P_0))\subseteq \mathcal{O}_{K, S}.\]
That is, $\mathfrak{a}_n$ is the largest ideal such that $f^n(P)\equiv P$ modulo ${\mathfrak{a}_n}$.
Then $\mathfrak{a}_n$ is a divisibility sequence, i.e., $m\mid n\Rightarrow \mathfrak{a}_m\mid\mathfrak{a}_n$, which we will call a \emph{H\'{e}non divisibility sequence}. In the case $\mathcal{O}_{K, S}=\mathbb{Z}$, we may identify the ideals with their unique positive generators, and think of this as a sequence of positive integers.

\begin{example} $f(x, y)=(y, y^2 + x - 2)$, $P=(2, 3)$.  
	\begin{multline*}
	\mathfrak{a}_n:1,
		1,
		8,
		3,
		1,
		8,
		1,
		3,
		8,
		1,
		5,
		48,
		11,
		1,
		8,
		51,
		1,
		8,
		1,
		3,
		8,
		5,
		7,
		288,
		13,
		11,
		8,
		3,
		1, ...	
	\end{multline*}
\end{example}

In analogy to other divisibility sequences, we probably expect this sequence to grow slowly. By comparing to the height of $P_n$, one easily obtains an upper bound of size $C^{d^n}$, for some $C$, on the norm of each of the terms in the gcd, but the gcd itself should usually be much smaller.
\begin{theorem}[Bugeaud, Corvaja, and Zannier~\cite{bcz}]
	If $a, b\geq 2$ are multiplicatively independent integers and $\varepsilon>0$, then
	\[\gcd(a^n-1, b^n-1)\ll_\varepsilon e^{\varepsilon n}.\]
\end{theorem}
\begin{theorem}[Huang~\cite{huang}]
	If $f(x), g(x)\in\mathbb{Z}[x]$ of degree $d\geq 2$, then ``generically'' and under Vojta's Conjecture
	\[\gcd(f^n(a)-\alpha, g^n(b)-\beta)\ll_\varepsilon e^{\varepsilon d^n}.\]
\end{theorem}
On the other hand, it is certainly true that every prime $\mathfrak{p}\subseteq\mathcal{O}_{K, S}$ divides some term in the sequence, since the image of $P$ in the residue field $\mathcal{O}_{K, S}/\mathfrak{p}$ must be periodic of period at most $\operatorname{Norm}(\mathfrak{p})^2$.

\begin{question} \normalfont
	Under what circumstance is it true that
	\[\log\operatorname{Norm}(\mathfrak{a}_n)=o(d^n) \]
	as $n\to \infty$? 
	Under what circumstance does there exist an ideal $\mathfrak{b} \subset \mathcal{O}_{K, S}$ such
that $\mathfrak{a}_n | \mathfrak{b}$ for infinitely many $n$? Note that, since $\mathfrak{a}_1 | \mathfrak{a}_n$ 
for all $n$, $\mathfrak{b} =\mathfrak{a}_1$ is a
reasonable candidate.
\end{question}

\bibliographystyle{alpha}
\bibliography{biblio1}
\end{document}